\documentclass[a4paper, 12pt]{article}

\usepackage[french, english]{babel}
\usepackage{amsmath}
\usepackage{a4}
\usepackage{euscript}
\usepackage{amsfonts}
\usepackage{amssymb}
\usepackage{xy}
\usepackage[dvips]{graphics}
\usepackage{enumerate}
\usepackage{array}
\usepackage{theorem}
\usepackage{epic}
\usepackage{color}
\usepackage{geometry}

\newtheorem{The}{Th\'eor\`eme}[subsection]
\newtheorem{Pro}[The]{Proposition}
\newtheorem{Def}[The]{D\'efinition}
\newtheorem{Lem}[The]{Lemme}
\newtheorem{Rem}[The]{Remarque}

\def\k#1{\kern#1em}
\def\Ib#1{{\rm I\kern-.25em#1}}
\def\Ibb#1{{\rm I\kern-.23em#1}}

\def\bs{{\baseskip}}
\def\ms{{\medskip}}
\def\sms{{\smallskip}}

\def\Eproof{{\hfill ///}}

\def\ZZ{{\mathbb Z}}
\def\QQ{{\mathbb Q}}

\def\PP{{\mathbb P}}

\def\lv{\mathbb }
\def\suml#1{{\textstyle\sum\limits_{#1}}}
\def\maxl#1{{\textstyle\max\limits_{#1}}}
\def\prodl#1{{\textstyle\prod\limits_{#1}}}

\def\mmu{\ensuremath{\boldsymbol \mu}}
\def\aalpha{\ensuremath{\boldsymbol \alpha}}

\newcommand{\Spec}{{\mathrm{Spec}\, }}
\newcommand{\Proj}{{\mathrm{Proj}\, }}

\newcommand{\rem}{{\mathrm{rem}}}
\newcommand{\Br}{{\mathrm{Br}}}
\newcommand{\Ram}{{\mathrm{Ram}\, }}
\newcommand{\rg}{{\mathrm{rang}}}

\newcommand{\disc}{{\mathrm{Disc}\, }}

\begin{document}

\title{Vers un algorithme pour la r\'eduction stable des rev\^etements $p$-cycliques de la droite projective sur un corps $p$-adique
\thanks{AMS classification 
number: 11G20, 14H30 (14Q05)}}
\author { Michel Matignon }
\maketitle
\begin{abstract}
Dans sa th\`ese C. Lehr vient d'exhiber un algorithme pour d\'ecrire la r\'eduction stable des 
rev\^etements $p$-cycliques de la droite projective sur un corps $p$-adique dans le cas o\`u le lieu
de branchement de cardinal $m+1$ a la g\'eom\'etrie \'equidistante et sous l'hypoth\`ese o\`u $m<p$. Dans cette note toujours dans le cas o\`u le lieu de branchement
a la g\'eom\'etrie \'equidistante nous proposons un algorithme sans condition sur $m$; en particulier nous 
pouvons \'etudier la r\'eduction en $2$ des courbes hyperelliptiques ayant un lieu de branchement \'equidistant.
\end{abstract}
\begin {abstract}
In his Ph. D. thesis, C. Lehr offers an algorithm which gives the stable model for $p$-cyclic covers of the projective line over a $p$-adic field under the conditions that the branch locus whose cardinal is $m+1$ has the so called equidistant  geometry and $m<p$. In this note we give an 
algorithm also in the equidistant geometry case but without condition on $m$. In particular we are able
to study the reduction at $2$ of hyperelliptic curves with equidistant branch locus.
\end{abstract}

\setcounter{section}{-1}
\section{Introduction}
\indent{\ \ \ }
Par corps $p$-adique nous entendons un corps $K$ d'in\'egale caract\'eristique $p>0$ qui est complet pour une valuation discr\`ete not\'ee $v$. On note $R$ l'anneau des entiers de $K$, $\pi$ une uniformisante et $k$ son corps r\'esiduel que l'on supposera par soucis de commodit\'e alg\'ebriquement clos. Nous supposerons de plus  que $R$ contient $\zeta$ une racine primitive $p$-i\`eme de l'unit\'e et noterons $\lambda=\zeta -1$. 

Les arithm\'eticiens savent bien qu'il est d\'elicat d'\'etudier la r\'eduction en $2$ d'une courbe 
elliptique ou plus g\'en\'eralement d'une courbe hyperelliptique. C'est l\`a un avatar d'un probl\`eme 
bien connu des g\'eom\`etres, \`a savoir le ``mauvais comportement'' en  $p$ du groupe fondamental 
alg\'ebrique.  Nous nous proposons d'examiner plus pr\'ecis\'ement la g\'eom\'etrie des rev\^etements $p$-cycliques de la droite projective sur $K$; par g\'eom\'etrie nous entendons l'\'etude du mod\`ele stable. 

Une premi\`ere tentative dans ce sens se trouve dans un expos\'e de s\'eminaire de Coleman 
([Co], \S 6); 
cette approche s'est d\'evelopp\'ee initialement en vue de comprendre la r\'eduction stable des
courbes de Fermat $F_n: X^n+Y^n=Z^n$ lorsque $p|n$; et avec succ\`es dans ce cas parce que l'on se ram\`ene \`a l'\'etude des rev\^etements cycliques de $\lv P^1$ ramifi\'es au-dessus de $3$ points. R\'ecemment il 
y a eu un regain d'int\'er\^et sur cette question principalement apr\`es les travaux de M. Raynaud 
qui a introduit dans [Ra 1] des m\'ethodes nouvelles  combinant th\'eorie de Galois, g\'eom\'etrie 
semi-stable et d\'eg\'en\'erescence des sch\'emas en groupes. Divers auteurs ont depuis contribu\'e \`a une description qualitative en termes combinatoires et diff\'erentiels de la 
g\'eom\'etrie $p$-adique de ces rev\^etements (cf. [Gr-Ma], [He], [Sa 1,2,3],...), la 
complexit\'e de la r\'eponse d\'ependant de la g\'eom\'etrie du lieu de branchement. Un retour sur 
l'aspect algorithmique a \'et\'e op\'er\'e tr\`es r\'ecemment dans sa th\`ese par C. Lehr. Soit  $C\to \lv P^1$ un rev\^etement $p$-cyclique de la droite projective 
ramifi\'e en $m+1\geq 3$ points. Il \'etudie le cas o\`u la g\'eom\'etrie du lieu de branchement est la 
plus simple d'un point de vue $p$-adique: les $m+1$ points sont \'equidistants i.e. pour un choix 
convenable de coordonn\'ees sur $\lv P^1$, ils sont dans des classes distinctes modulo $\pi$. Pour cette
g\'eom\'etrie et sous les  conditions suppl\'ementaires $(m,p)=1$ et $(m-1,p)=1$, C. Lehr donne un crit\`ere de potentielle bonne r\'eduction et, sous la condition $m<p$, un algorithme qui d\'ecrit le mod\`ele stable. Ces conditions sont tr\`es restrictives pour $p$ petit et en particulier si $p=2$, elles
\'eliminent toutes les courbes hyperelliptiques. 

Dans cette note, toujours sous l'hypoth\`ese que le lieu de branchement a la g\'eom\'etrie 
\'equidistante, nous proposons un algorithme qui fonctionne sans condition. 

Passons \`a une description de la note. Nous esp\'erons que le lecteur partagera le choix de la terminologie et des  notations. 

Dans la partie 1, nous rappelons les propri\'et\'es essentielles \`a notre propos concernant la d\'eg\'en\'ere-scence des  $\mmu_p$-torseurs et la g\'eom\'etrie semi-stable. Nous pr\'esentons bri\`evement les algorithmes de Coleman et de Lehr afin de pr\'eparer le lecteur \`a l'introduction dans la 
partie 2  d'une g\'en\'eralisation des d\'eveloppements de Taylor des 
polyn\^omes sur un corps $p$-adique. Bri\`evement, pour un entier $m$ donn\'e, pour $F(X)\in R[X]$ et 
$y\in R^{alg}$, on d\'eveloppe $F(X+y)\in R[y][X]$ sous la forme d'une approximation dans $R^{alg}[X]$ \`a presque 
$\frac{p}{p-1}v(p)$ pr\`es par la puissance $p$-i\`eme d'un polyn\^ome de mani\`ere que les coefficients des mon\^omes restant de degr\'e $\leq m$ et multiple de $p$ aient  une taille n\'egligeable dans le sens de l'approximation. Le coefficient $F^{[1]}(y)$ du mon\^ome $X$ est une ``fonction multiforme de $y$'' dont la  norme  
$N_m(F^{[1]}(y))$ est une fonction polyn\^ome de $y$ qui modulo $p$ est une puissance $F'(y)^{p^{r(m)}}$ de la d\'eriv\'ee usuelle, pour $r(m)$ convenable. 

Soit $\Br$, le lieu de branchement du rev\^etement suppos\'e de cardinal $m+1\geq 3$, et  $Z^p=F(X)\in R[X]$ une \'equation du $\mmu_p$-torseur  au-dessus de $\lv P^1-\Br$. Dans ([Le 2], Th. 4.1) Lehr a montr\'e, sous l'hypoth\`ese $m<p$, que certains  z\'eros de $F'(Y)$ donnent des centres 
pour les disques ferm\'es de $\lv P^1_K$ qui induisent des composantes de genre non nul dans le 
mod\`ele stable du rev\^etement. La preuve n'est pas difficile et vient d'une comparaison des 
polygones de Newton de $F(X+y)$  et $F'(X+y)$. Dans le cas g\'en\'eral $F'(Y)$ est remplac\'e 
par $N_m(F^{[1]}(Y))$ et l'analyse des z\'eros se fait gr\^ace \`a une \'etude fine de la 
d\'eg\'en\'erescence du $\mmu_p$-torseur; nous renvoyons le lecteur au th\'eor\`eme 3.2.2 pour un \'enonc\'e 
pr\'ecis. 

La partie 4 est consacr\'ee aux exemples et en premier lieu \`a la situation o\`u les points de 
branchement sont en position \'equidistante et leurs classes modulo $\pi$ en position g\'en\'erale.
Puis le cas des courbes hyperelliptiques ($p=2$) est plus pr\'ecis\'ement \'etudi\'e pour $m$ petit.

La partie 5 donne une  g\'en\'eralisation dans le cas des rev\^etements $p$-cycliques d'une courbe ayant bonne r\'eduction. 

\section{G\'en\'eralit\'es}
\subsection{R\'eduction des $\mmu_p$-torseurs et la diff\'erente}
\indent{\ \ \ }
On reprend les notations de l'introductions. 
La proposition qui suit est un r\'esum\'e commode pour les applications  que nous avons en vue; elle 
se retrouve dans les travaux de divers auteurs; nous renvoyons  \`a ([He], ou [Sa 3]) pour un expos\'e complet. 
Nous devons rappeler la d\'efinition du sch\'ema en groupe ${\cal H}_n$. 

Pour tout entier $n>0$, on note ${\cal G}_n:=\Spec R[X,\frac{1}{\pi^nX+1}]$ dont la fibre g\'en\'erique 
est isomorphe au groupe multiplicatif et la fibre sp\'eciale 
s'identifie au groupe additif.  
Pour $0<n\leq v_K(\lambda)$, le polyn\^ome $\frac{(\pi^nX+1)^p-1}{\pi^{pn}}$ est \`a coefficients
dans $R$; l'homomorphisme $$\Psi_n: R[Y,\frac{1}{\pi^{pn}Y+1}]\to R[X,\frac{1}{\pi^nX+1}],$$ d\'efini par $\Psi_n(Y)=\frac{(\pi^nX+1)^p-1}{\pi^{pn}}$  est une isog\'enie de degr\'e $p$; on note ${\cal H}_n$ le 
noyau de $\Psi_n$. Le sch\'ema  ${\cal H}_n$ est fini plat sur $R$, de degr\'e $p$. Sa fibre g\'en\'erique est isomorphe au groupe $\mmu_{p,K}$. Si $0<n<v_K(\lambda)$, sa fibre sp\'eciale est le groupe 
radiciel additif $\aalpha_p$ et si $n=v_K(\lambda)$, sa fibre sp\'eciale est le groupe \'etale isomorphe 
\`a $\ZZ/p\ZZ$. 

\begin{Pro} 
\label{StructMup} 
{\bf(\cite{He}, prop. 1.6)}.--- Soit $X:=\Spec A$ un sch\'ema affine plat sur $R$, dont
les fibres sont int\`egres et de dimension $1$; on suppose que $A$ est une $R$-alg\`ebre factorielle 
et compl\`ete pour la topologie $\pi$-adique. Soit $Y_K\to X_K$ un $\mmu_p$-torseur \'etale non trivial,
donn\'e par une \'equation $y^p=f$, o\`u $f$ est inversible dans $A_K$, et $Y$ le normalis\'e de $X$
dans $Y_K$; on suppose que la fibre sp\'eciale de $Y$ est int\`egre. Soit $\eta$ (resp. $\eta'$) le point g\'en\'erique de la fibre sp\'eciale de $X$ (resp. $Y$). Les anneaux locaux ${\cal O}_{X,\eta}$ et
${\cal O}_{Y,\eta'}$ sont alors des anneaux de valuation discr\`ete d'uniformisante $\pi$. Notons
$\delta$ la valuation de la diff\'erente de ${\cal O}_{Y,\eta'}/{\cal O}_{X,\eta}$. On distingue alors
deux cas suivant la valeur de $\delta$.
\begin{enumerate}[-]
\item 
Si $\delta=v_K(p)$, $Y$ est un $\mmu_{p,R}$-torseur pour la topologie fppf,
donc $Y=\Spec B$, avec $B:=\frac{A[Y]}{(Y^p-u)}$, o\`u $u$ est une unit\'e de $A$, unique \`a la multiplication d'une puissance $p$-i\`eme d'une unit\'e de $A$ pr\`es. On dit que le torseur a {\bf r\'eduction multiplicative}.
\item Si $0\leq \delta <v_K(p)$, on a $\delta=v_K(p)-n(p-1)$, o\`u $n$ est un entier tel que 
$0\leq \delta <v_K(\lambda)$, et $Y\to X$ est un torseur sous ${\cal H}_n$ pour la topologie fppf, donc donn\'e par $B:=\frac{A[W]}{\frac{(\pi^nW+1)^p-1}{\pi^{pn}}-u},$ o\`u $u$ est un \'el\'ement de 
$A$. De plus, si $B$ est isomorphe \`a $\frac{A[W]}{\frac{(\pi^nW+1)^p-1}{\pi^{pn}}-u'}$, il existe
$v\in A$ tel que $u'=u(\pi^nv+1)^p+\frac{(\pi^nv+1)^p-1}{\pi^{pn}}.$ 
\end{enumerate}
Si $0\leq \delta <v_K(p)$ (resp. $\delta=0$), on dit que le torseur $Y\to X$ a {\bf r\'eduction 
additive} (resp. {\bf r\'eduction \'etale}). 
\end{Pro}

\subsection{Mod\`ele stable} 
\indent{\ \ \ }
Nous adoptons la terminologie de [Li] \`a savoir qu'une $R$-courbe propre est {\it semi-stable} (resp. {\it stable}) si ses fibres 
g\'eom\'etriques sont des courbes semi-stables (resp. stables), i.e. sont projectives r\'eduites et connexes et  ont pour seules singularit\'es des points doubles ordinaires (resp. semi-stables, de genre arithm\'etique $\geq 2$ et  chacune de leurs composantes 
irr\'eductibles isomorphe \`a $\lv P^1$ rencontre les autres composantes en au moins 3 points). 

Une $R$-courbe propre est appel\'ee mod\`ele de sa fibre g\'en\'erique. 

On \'etend la notion de mod\`ele  stable  au contexte des courbes lisses propres 
sur $K$  
point\'ees par un ensemble fini $S\subset C(K)$. On appelle {\it mod\`ele stable de la courbe point\'ee 
$(C,S)$} un mod\`ele ${\cal C}/R$ de $C$ tel que  les points de $S$ se sp\'ecialisent en des points distincts $\bar S$ (i.e. la cl\^oture sch\'ematique de  $S$ est lisse sur $R$) du lieu
lisses de  ${\cal C}_s$ et tel que  les fibres sont 
semi-stables  et chacune de leurs composantes 
irr\'eductibles isomorphe \`a $\lv P^1$ contient au moins 3 points parmi les points 
de $\bar S$ ou les points d'intersection.  

Il suit du th\'eor\`eme de r\'eduction semi-stable que 
si $C$, resp. $(C,S)$ est une courbe (resp. point\'ee) propre, lisse et g\'eom\'etriquement irr\'eductible sur $K$, de genre $\geq 2$ (resp. de genre $g$ avec $2g+|S|-1\geq 2$) 
il existe une extension finie $K'/K$ telle que $C_{K'}:=C\times _KK'$ admet un mod\`ele stable sur 
la cl\^oture int\'egrale $R'$ de $R$ dans $K'$ (voir [Ab] ou [Li],  pour une d\'emonstration du th\'eor\`eme de r\'eduction
semi-stable). 

Dans ce qui suit $(C,G)$ d\'esigne  une courbe $C/K$ munie d'une action par un groupe de $K$-automorphismes 
$G$; on note  $\Ram$  (resp. $\Br$) le lieu de ramification (resp. branchement) de $f: C\to C/G$; nous dirons que $(C,G)$ admet un  mod\`ele stable sur $R$ si $\Ram\subset C(K)$ et si la courbe point\'ee 
$(C,\Ram)$ admet un mod\`ele stable sur $R$; nous appelerons alors {\bf mod\`ele stable de $(C,G)$ le mod\`ele stable de la courbe point\'ee 
$(C,\Ram)$}. Supposons pour simplifier  que ce mod\`ele stable est d\'efini sur $R$; alors  par 
l'unicit\'e il est \'equivariant et le morphisme $f:C\to C/G$ s'\'etend en un morphisme fini de 
$R$-courbes  $ f:  {\cal C}\to {\cal C}/G$ o\`u ${\cal C}/G$ est semi-stable et s'obtient \`a partir 
du mod\`ele stable de $ C/G$ point\'ee par $\Br$ en \'eclatant des points ferm\'es de la fibre 
sp\'eciale; en particulier sa fibre sp\'eciale ne diff\`ere de celle du mod\`ele stable que par 
des $\lv P^1_k$ ([Ra 1]). 
\begin{Def}.--- {\rm 
Nous dirons que le lieu de branchement $\Br$ a la {\it g\'eom\'etrie \'equidistante} si le mod\'ele 
stable de la courbe point\'ee $(C/G,\Br)$ est d\'efini et lisse sur $R$ i.e. $C/G$ admet un mod\`ele 
lisse sur $R$ et les points de $\Br$ se sp\'ecialisent dans ce mod\`ele en des 
points distincts. }
\end{Def}
\subsection{L'algorithme de Coleman} 
\indent{\ \ \ }
Soit $C$ une courbe propre et lisse sur $K$. Soit $n>2$ un entier diff\'erent de $p$, on sait (c.f. [Ab], prop. 5.10) que $C$ a un mod\`ele stable sur l'extension $K_{n}$ de $K$ obtenue en adjoignant les points de $n$-torsion de sa jacobienne pour un $n$ entier diff\'erent de $p$. 
Malheureusement \'etant donn\'ee une
\'equation de $C$ il n'est pas facile de trouver $K_n$ et donc de trouver le mod\`ele stable de $C$. Dans ([Co], \S  6) Coleman dit: ``It would be desirable to have some efficient algorithm for computing
the stable reduction of a curve.'' 

Coleman est principalement interess\'e  par les rev\^etements galoisien de $\lv P^1$ et constate que 
le cas difficile est celui o\`u $p$ divise l'ordre du groupe. Il d\'eveloppe un algorithme dans
le cas $p$-cyclique qui marche dans certains cas. Pr\'ecis\'ement, il consid\`ere une \'equation 
$Z^p=F(X)\in K[X]$ avec $\deg F=n$ et $V/K$ la vari\'et\'e  affine sur $K$, 
$\Spec \frac{K[Y,a_i,b_j,\, 0\leq i\leq [\frac{n}{p}], 0<j\leq n,\, (j,p)=1]}{\cal I}$, o\`u l'id\'eal $\cal I$  des \'equations est cod\'e par l'identit\'e
$$(*)\  F(X+Y)=(\suml{0\leq i\leq [\frac{n}{p}]}a_iX^i)^p+\suml{0<i\leq n,\, (i,p)=1}b_iX^i.$$
M\^eme si ce syst\`eme a trivialement une solution modulo $p$, il n'est pas s\^ur qu'il a des solutions; en effet un d\'ecompte montre que la dimension de $V$ sur $K$ est en g\'en\'eral \'egale \`a $1$. 
Si l'on impose la condition suppl\'ementaire $b_1=0$; on obtient un ensemble \'eventuellement vide et 
au plus fini de valeurs de $Y$ pour lesquelles l'\'equation (*) permet parfois d'exhiber un mod\`ele entier de la courbe $Z^p=F(X)$ qui est propice \`a donner des composantes de genre non nul en r\'eduction; ces composantes induisent des disques ferm\'es de la droite projective quotient et les z\'eros de $b_1$   en fournissent des centres. La recherche 
de centres de ces disques est un probl\`eme de nature analytique  et la m\'ethode de Coleman en 
propose une attaque alg\'ebrique (recherche d'une vari\'et\'e de dimension $0$ et \'eventuellement vide) aussi cette  m\'ethode doit \^etre modifi\'ee pour pouvoir esp\'erer aboutir dans le cas g\'en\'eral. Nous renvoyons le lecteur \`a 4.2 pour des exemples trait\'es par la m\'ethode de Coleman. 
\subsection{L'algorithme de Lehr}
\indent{\ \ \ }
Dans sa th\`ese C. Lehr a \'etudi\'e la r\'eduction stable des rev\^etements $p$-cycliques de 
la droite projective sur un corps $p$-adique dans le cas o\`u le lieu de branchement a la g\'eom\'etrie \'equidistante.  

Soit $Z^p=F(X)$ une \'equation normalis\'ee du rev\^etement $C\to \lv P^1$ i.e.
 $F(X)=\prodl{1\leq i\leq m}
(X-x_i)^{e_i}\in R[X]$ de degr\'e $N$ avec $(N,p)=1$  et $(e_i,p)=1$ pour $1\leq i\leq m$. On suppose 
de plus que $v(x_i)=v(x_i-x_j)=0$ pour tout $i\neq j$. Ainsi le 
lieu de branchement  $\Br =\{\infty ,x_1,..,x_m\}$  a la g\'eom\'etrie \'equidistante. Dans le cas 
o\`u $m<p$, C. Lehr montre que le d\'eveloppement de Taylor $F(X+y)=F(y)+F'(y)X+...$ pour 
$y\in R^{alg}$ un z\'ero de $F'(Y)$ avec $v(F(y))=0$ donne naissance aux composantes du mod\`ele 
stable de $C$ qui sont de genre non nul (cf. [Le 2], Th. 4.1).

\section{$p$-d\'eveloppements de Taylor}
\indent{\ \ \ }
Les champs respectifs de succ\`es des algorithmes de 
Coleman et Lehr nous am\`enent \`a g\'en\'eraliser la notion de d\'eveloppement de Taylor qui approche dans un sens $p$-adique le d\'eveloppement id\'eal de Coleman et pr\'esente  l'avantage de ne pas \^etre 
sensible au petites variations des coefficients.
\subsection{D\'efinition}
\begin{Def}.--- {\rm On fixe $K(Y)^{alg}$ une cl\^oture alg\'ebrique de $K(Y)$ munie d'une valuation $v_Y$ qui prolonge la valuation de Gauss de $K(Y)$ relative \`a $Y$ et $R[Y]^{alg}$ d\'esigne la cl\^oture int\'egrale de $R[Y]$. Soit $m,n$ deux entiers; on note $c_n:=1+1/p+...+1/p^n$
et $I_{m,n}$ l'id\'eal $(p^{c_n}X,X^{m+1})$ de $R[Y]^{alg}[X]$. 
Un {\it $p$-d\'eveloppement de Taylor d'ordre $m$ et de niveau $n$ de 
$F(X)\in R[X]$}  est la donn\'ee d'un couple $(A_{m,n}(X,Y),B_{m,n}(X,Y))$ de polyn\^omes $\in R[Y]^{alg}[X]$ avec
$$A_{m,n}(X,Y)=\suml{0\leq i\leq m/p}a_i(Y)X^i, \ B_{m,n}(X,Y)=\suml{1\leq j\leq m,\ (j,p)=1}b_j(Y)X^j$$  
tels que 
$$F(X+Y)-(A_{m,n}(X,Y))^p-B_{m,n}(X,Y)\in I_{m,n}.$$}
\end{Def}
 
\begin{Rem} .---\begin{enumerate} [i)]{\rm
\item Si $m<p$ le d\'eveloppement de Taylor classique induit un $p$- d\'eveloppement de Taylor d'ordre $m$ et de niveau $n$ pour tout $n$.
\item Notez que $c_nv(p)=(1-1/p^{n+1})v(\lambda^p)$; ainsi la formule pr\'ec\'edente 
donne modulo $X^{m+1}$ une approximation de $F(X+Y)$ par une puissance $p$-i\`eme \`a presque  $v(\lambda ^p)=\frac{p}{p-1}v(p)$ pr\`es.}
\end{enumerate}
\end{Rem}
\subsection{Existence}
Si $L$ est un corps alg\'ebriquement clos et si $x\in L$, on notera $x^{1/p}$ une racine 
$p$-i\`eme de $x$.
On a la proposition suivante:
\begin{Pro} .--- On fixe l'entier $m$, notons $F_m(X+Y)=s_0(Y)+s_1(Y)X+....+s_m(Y)X^m\in R[Y][X]$ le d\'eveloppement de Taylor \`a l'ordre $m$ de $F(X)\in R[X]$. On d\'efinit une  suite de 
polyn\^omes $ \in R[Y]^{alg}[X]$, $(A_{m,n}(X,Y),B_{m,n}(X,Y))$ par les relations de 
r\'ecurrence suivantes:
$$(  \rg \  0)\quad  (A_0(X,Y),B_0(X,Y))=(\suml{0\leq i\leq m/p}s_{ip}(Y)^{1/p}X^i,\suml{1\leq j\leq m,\ (j,p)=1}s_j(Y)X^j).$$
Pour d\'efinir la suite au rang $(n+1)$, on consid\`ere le reste au rang $n$, $R_n(X,Y):=F_m(X+Y)-(A_{m,n}(X,Y))^p-B_{m,n}(X,Y)=p^{c_n}(r_1(Y)X+...+r_m(Y)X^m)$ o\`u $r_i(Y)\in R[Y]^{alg}$. Soit 
$AR_n:=p^{c_{n}/p}\suml{1\leq i\leq m/p}r_{ip}(Y)^{1/p}X^i$ et $BR_n:=p^{c_n}\suml{1\leq j\leq m,\ (j,p)=1}r_j(Y)X^j$, alors:
$$(A_{m,n+1}(X,Y),B_{m,n+1}(X,Y)):=(A_{m,n}(X,Y)+AR_n,B_{m,n}(X,Y)+BR_n).$$
Les polyn\^omes $(A_{m,n}(X,Y),B_{m,n}(X,Y))$ d\'efinissent un $p$-d\'eveloppement de Taylor d'ordre $m$ et de niveau $n$ de $F(X).$
\end{Pro}
{\it Preuve.} Avant de passer \`a la preuve remarquons que la suite pr\'ec\'edente d\'epend \`a chaque 
\'etape du choix de racines $p$-i\`emes et l'\'enonc\'e vaut donc pour un choix quelconque. L'\'etape $0$
est une \'egalit\'e en caract\'eristique $p$. Les polyn\^omes mis en jeu \`a chaque \'etape sont \`a coefficients dans l'anneau $ R[Y]^{alg}$ muni de la valuation de Gauss $v_Y$. Il reste \`a voir qu'au rang $n+1$ on a une approximation 
modulo $p^{c_{n+1}}$. On calcule le reste $R_{n+1}(X,Y)=F_m(X+Y)-(A_{m,n+1}(X,Y))^p-B_{m,n+1}(X,Y)=
F_m(X+Y)-(A_{m,n}(X,Y)+AR_n)^p-B_{m,n}(X,Y)-BR_n=S_1+S_2$
o\`u $S_1:=R_n(X,Y)-AR_n^p-BR_n$ et $S_2=A_{m,n}(X,Y)^p+AR_n^p-(A_{m,n}(X,Y)+AR_n)^p$. Puisque $v_Y(R_n(X,Y))\geq c_nv(p)$ la divisibilit\'e par $p$ des coefficients binomiaux implique que $v_Y(A)\geq (c_n+1)v(p)$. De m\^eme $v_Y(B)\geq v(p)+max(v_Y(A_{m,n}),v_Y(AR_n))\geq (1+c_n/p)v(p)=c_{n+1}v(p).$

\Eproof

\begin{Rem}.---  {\rm Puisque les coefficients des polyn\^omes $(A_{m,n},B_{m,n})$ sont obtenus par 
extraction de racines $p$-i\`emes, ils vivent dans  
une extension  $L_n/K(Y)$  galoisienne  de groupe $G_n$, un $p$-groupe. Si $\tau \in G_n$, on note  
$(A_{m,n}^{\tau},B_{m,n}^{\tau})$ les polyn\^omes obtenus par l'action de $\tau$ 
sur les coefficients; ils donnent donc naissance \`a un $p$-d\'eveloppement de Taylor d'ordre $m$ et de niveau $n$ de $F(X).$ }
\end{Rem}
\subsection{Unicit\'e}
\indent{\ \ \ }
Afin de comparer les $p$-d\'eveloppements de Taylor d'un polyn\^ome nous \'etablissons un lemme:
\begin{Lem} .--- Soit $\cal R$, un anneau local dominant $R^{alg}$. Soient $a_i,a'_i,b_j,b'_j\in 
\cal R $ alors $$(\suml{1\leq i\leq [m/p]}a_iX^i)^p+\suml{1\leq j\leq m,\ (j,p)=1}b_jX^j=
(\suml{1\leq i\leq [m/p]}a'_iX^i)^p+\suml{1\leq j\leq m,\ (j,p)=1}b'_jX^j \,  \mod (p^{c_n})$$
si et seulement si  $a'_i=a_i \, \mod (p^{c_{n}/p})$ et $b'_j=b_j \,  \mod (p^{c_n})$.
\end{Lem}
  
{\it Preuve.} Supposons que  $a'_i=a_i+p^{c_{n}/p}a''_i$ et $b'_j=b_j +p^{c_n}b''_j$, on a 
$$(\suml{1\leq i\leq [m/p]}a'_iX^i)^p+\suml{1\leq j\leq m,\ (j,p)=1}b'_jX^j =(\suml{1\leq i\leq [m/p]}(a_i+p^{c_{n}/p}a''_i)X^i)^p+\suml{1\leq j\leq m,\ (j,p)=1}(b_j+p^{c_n}b''_j)X^j=$$
$$=(\suml{1\leq i\leq [m/p]}a_iX^i+p^{c_{n}/p}(\suml{1\leq i\leq [m/p]}a''_iX^i))^p
+\suml{1\leq j\leq m,\ (j,p)=1}b_jX^j \mod (p^{c_n})=$$
$$=(\suml{1\leq i\leq [m/p]}a_iX^i)^p+\suml{1\leq j\leq m,\ (j,p)=1}b_jX^j \mod (p^{c_n})$$
puisque $pp^{c_{n}/p}=p^{c_{n+1}}\in  (p^{c_n})$.

L'implication r\'eciproque se fait par r\'ecurrence sur $n$. Si $n=0$, on a une \'egalit\'e modulo $p$.
Supposons donc avoir une \'egalit\'e 

$$(\suml{1\leq i\leq [m/p]}a_iX^i)^p+\suml{1\leq j\leq m,\ (j,p)=1}b_jX^j=
(\suml{1\leq i\leq [m/p]}a'_iX^i)^p+\suml{1\leq j\leq m,\ (j,p)=1}b'_jX^j \,  \mod (p^{c_{n+1}})$$
la r\'ecurrence permet d'\'ecrire  $a'_i=a_i+p^{c_{n}/p}a''_i$ et $b'_j=b_j +p^{c_n}b''_j$;
et le calcul pr\'ec\'edent donne la congruence 
$$(\suml{1\leq i\leq [m/p]}a'_iX^i)^p+\suml{1\leq j\leq m,\ (j,p)=1}b'_jX^j =$$
$$=(\suml{1\leq i\leq [m/p]}a_iX^i)^p+p^{c_{n}}(\suml{1\leq i\leq [m/p]}a''_iX^i)^p+ \suml{1\leq j\leq m,\ (j,p)=1}b_jX^j +p^{c_n}\suml{1\leq j\leq m,\ (j,p)=1}b''_jX^j \,  \mod (p^{c_{n+1}})$$
et donc 
$$p^{c_{n}}(\suml{1\leq i\leq [m/p]}a''_iX^i)^p+p^{c_n}\suml{1\leq j\leq m,\ (j,p)=1}b''_jX^j =0\,  \mod (p^{c_{n+1}})$$ 
autrement dit 
$$(\suml{1\leq i\leq [m/p]}a''_iX^i)^p+\suml{1\leq j\leq m,\ (j,p)=1}b''_jX^j =0\,  \mod 
(p^{1/p^{n+1}})$$ et donc
$$\suml{1\leq i\leq [m/p]}{a''_i}^{p}X^{ip}+\suml{1\leq j\leq m,\ (j,p)=1}b''_jX^j =0\,  \mod 
(p^{1/p^{n+1}})$$ ainsi $a''_i=0\,  \mod (p^{1/p^{n+2}})$ et $b''_j=0\,  \mod (p^{1/p^{n+1}}).$
\Eproof

\begin{Def} .--- {\rm Soit $m,n$ deux entiers. On dit que 2 couples de polyn\^omes $\in R[Y]^{alg}[X]$
$$(A_{m,n}(X,Y)=\suml{0\leq i\leq m/p}a_i(Y)X^i, \ B_{m,n}(X,Y)=\suml{1\leq j\leq m,\ (j,p)=1}b_j(Y)X^j),$$
$$(A'_{m,n}(X,Y)=\suml{0\leq i\leq m/p}a'_i(Y)X^i, \ B'_{m,n}(X,Y)=\suml{1\leq j\leq m,\ (j,p)=1}b'_j(Y)X^j),$$
sont \'equivalents si et seulement si $$(A_{m,n}(X,Y))^p+B_{m,n}(X,Y)=(A'_{m,n}(X,Y))^p+B'_{m,n}(X,Y) \mod (p^{c_n}).$$ }
\end{Def}
\begin{Rem} .--- {\rm Revenons \`a la remarque 2.2.2, les couples $(A_{m,n},B_{m,n})$ et 
$(A_{m,n}^{\tau},B_{m,n}^{\tau})$ sont donc \'equivalents. On peut montrer en revenant \`a la 
construction que $v_Y( A_{m,n} -A_{m,n}^{\tau})\geq \frac{c_n}{p}v(p)+\frac{1}{p}v(\lambda)$ et 
$v_Y( B_{m,n} -B_{m,n}^{\tau})\geq c_nv(p)+v(\lambda).$}
\end{Rem}

 Soit $F(X)\in R[X]$; par le lemme 2.3.1, deux couples de polyn\^omes $\in R[Y]^{alg}[X]$
d\'efinissant un $p$-d\'eveloppement de Taylor d'ordre $m$ et de niveau $n$ de  $F(X)\in R[X]$ sont donc
\'equivalents; nous sommes ainsi amen\'es \`a la 
\begin{Def} .--- {\rm Soit $m,n$ deux entiers, on appelle {\it $p$-d\'eveloppement de Taylor sp\'ecial d'ordre $m$ et de niveau $n$ de  $F(X)\in R[X]$}
un couple $(A_{m,n},B_{m,n})$ construit comme dans la proposition 2.2.1, en particulier  le terme 
constant dans $A_{m,n}$ vaut $F(Y)^{1/p}$. Par analogie avec le d\'eveloppement de Taylor classique 
nous \'ecrirons 
$$(*) F(X+Y)=E(X,Y)^p+\suml{1\leq j\leq m,\ (j,p)=1}F^{[j]}(Y)X^j \mod (p^{c_n}X,X^{m+1})$$
o\`u $E(X,Y)\in (R[Y])^{alg}[X]$,  $\deg_X E(X,Y)\leq [m/p]$ et $F^{[j]}(Y)\in (R[Y])^{alg}$ (notez que la congruence impose l'\'egalit\'e $E(0,Y)^p=F(Y)$).

On appelle {\it groupe de Galois attach\'e au $p$-d\'eveloppement de Taylor (*)} le groupe de Galois
 $G_n$ de la $p$-extension $L_n$ de $K(Y)$ obtenue en adjoignant les racines $p$-i\`emes n\'ecessaires 
au $p$-d\'eveloppement de Taylor (cf. remarque 2.2.2).} 
\end{Def}
\subsection
{Sp\'ecialisation}
\indent{\ \ \ }
Nous serons amen\'es \`a sp\'ecialiser la formule $(*)$; pr\'ecis\'ement si $y\in R^{alg}$ on peut 
tester sur $F(X+y)$ l'algorithme de la proposition 2.2.1, l'interp\'etation que l'on doit donner 
au d\'eveloppement ainsi construit est la suivante:  quitte \`a faire une extension finie on peut supposer 
que $y\in R$; soit $\cal P$ une place au-dessus de $(Y-y)$ du corps de fonctions $L_n/K$ , on
peut regarder l'image de la formule $(*)$  dans le corps r\'esiduel $K({\cal P})$; on a  
une congruence 

$$(**) F(X+y)=E(X,Y)^p({\cal P})+\suml{1\leq j\leq m,\ (j,p)=1}F^{[j]}(Y)({\cal P})X^j \mod (p^{c_n}X,X^{m+1})$$ 
\noindent
que par abus de notation compatible avec la stabilit\'e par conjugaison nous \'ecrirons   
$$(***) F(X+y)=E(X,y)^p+\suml{1\leq j\leq m,\ (j,p)=1}F^{[j]}(y)X^j \mod (p^{c_n}X,X^{m+1}).$$

Les z\'eros de $F^{[1]}(Y)\in (R[Y])^{alg}[X]$ vont jouer un r\^ole d\'eterminant dans l'\'etude 
de la r\'eduction stable des torseurs $Z^p=F(X)$. Puisque $F^{[1]}(Y)$ n'est bien d\'efini que modulo l'action du groupe $G_n$ 
(cf. remarque 2.2.2) nous introduisons la norme $N_{L_n/K(Y)}(F^{[1]}(Y))\in R[Y]$ qui modulo $p$ est une puissance $|G_n|$-i\`eme de la d\'eriv\'ee $F'(Y)$; en fait, $m$ sera d\'etermin\'e par le lieu de 
branchement du torseur et $n$ sera l'entier tel que $p^n\leq m<p^{n+1}$. Ainsi on est amen\'e \`a d\'efinir:
\begin{Def} .--- {\rm Soit $m$ un entier naturel  et $n$,  l'entier tel que $p^n\leq m<p^{n+1}.$  
Nous appellerons {\it $p$-d\'eveloppement de Taylor sp\'ecial d'ordre $m$} de $F(X)\in R[X]$ un 
$p$-d\'eveloppement de Taylor sp\'ecial d'ordre $m$ et de niveau $n$
$$(*) F(X+Y)=E(X,Y)^p+\suml{1\leq j\leq m,\ (j,p)=1}F^{[j]}(Y)X^j \mod (p^{c_n}X,X^{m+1})$$
 et nous appellerons {\it $p$-d\'eriv\'ee de niveau $m$}  de $F(X)\in R[X]$, le polyn\^ome  
 $$N_m(F^{[1]}(Y)):=N_{L_n/K(Y)}(F^{[1]}(Y))\in R[Y].$$} 
\end{Def}
\section
{$p$-d\'eveloppements de Taylor et r\'eduction stable}
\indent{\ \ \ }
Dans cette partie on fixe un rev\^etement $p$-cyclique $C\to \lv P^1$ de la droite projective sur $K$
ramifi\'e en au moins 3 points (i.e. $g(C)>0$). 
Nous nous proposons de d\'ecrire le mod\`ele stable de $(C,G=\ZZ/p\ZZ)$ dans le cas o\`u le lieu de branchement a la g\'eom\'etrie \'equidistante (cf. 1.2.1). Fixons les notations.
\subsection
{Notations}
\begin{enumerate}[i)]\item Si $P(X)\in R[X]$ on note $\bar P[X]$ son image dans $k[X]$.
\item Soit $Z^p=F(X)$ une \'equation normalis\'ee du rev\^etement $C\to \lv P^1$ i.e. 
 $F(X)=\prodl{1\leq i\leq m}
(X-x_i)^{e_i}\in R[X]$ de degr\'e $N$ avec $(N,p)=1$  et $(e_i,p)=1$ pour $1\leq i\leq m$. On suppose 
de plus que $x_i\in R$ et que $v(x_i)=v(x_i-x_j)=0$ pour tout $i\neq j$. Ainsi le 
lieu de branchement $\Br =\{\infty ,x_1,..,x_m\}$  a la g\'eom\'etrie \'equidistante. On note 
$\bar \Br =\{\infty ,\bar x_1,..,\bar x_m\}$ la sp\'ecialisation de $\Br$. 
\item Soit $n$,  l'entier tel que $p^n\leq m<p^{n+1} $ et $c_n:=1+1/p+...+1/p^n$, enfin $N_m(F^{[1]}(Y))=\prodl{\tau \in G_n}{F^{[1]}}^{\tau}(Y)$ d\'esigne  la $p$-d\'eriv\'ee de niveau $m$ de $F(X)$ (cf. 2.4.1)  o\`u $G_n$ est d\'efini dans 2.3.4 et $|G_n|=p^{r(m)}$.
 \end{enumerate}
\begin{Pro} .--- La $p$-d\'eriv\'ee de niveau $m$ de $F(X)$, $N_m(F^{[1]}(Y))$ appartient \`a $R[Y]$;  de plus 
$N_m(F^{[1]}(Y))=\bar F'(Y)^{p^{r(m)}} \mod p$ et $\deg N_m(F^{[1]}(Y))=\deg \bar F'(Y)^{p^{r(m)}}$. 
\end{Pro}
{\it Preuve.} Montrons l'in\'egalit\'e $\deg N_m(F^{[1]}(Y))\leq \deg \bar F'(Y)^{p^{r(m)}}$. L'id\'ee est d'introduire un poids dans $R[Y]^{alg}$ qui g\'en\'eralise le degr\'e dans $R[Y]$ et de suivre dans l'algorithme qui donne le $p$-d\'eveloppement de Taylor sp\'ecial d'ordre $m$ et de niveau $n$ le poids du coefficient $b_1(Y)$ de $X$  dans  $B_{m,n}(X,Y)$. 

Pr\'ecis\'ement $-\deg_Y$ est la valuation discr\`ete de $K(Y)$ au point $\infty$; cette valuation 
s'\'etend  sur $K(Y)^{alg}$ en une semi-norme $w$: si $f\in K(Y)^{alg}$, soit 
$P_f(X):=\suml{0\leq i\leq N} a_iX^i$ le polyn\^ome irr\'eductible unitaire de $f$ sur $K(Y)$; alors $w(f):=\inf_i \frac {1}{N-i}(-\deg_Y a_i)$. 
On reprend les notations de la prop. 2.2.1. On a $F_m(X+Y)=s_0(Y)+s_1(Y)X+...$; on remarque 
que $w(s_i(Y))\geq (i-N)$ o\`u $N=\deg_X F$ est premier \`a $p$, ainsi $i-w(s_i(Y))\leq N$. Consid\'erons les formules au rang $0$. 
On a $w((s_{ip}(Y))^{1/p})=(ip-N)/p=i-N/p$ et donc  $i-w((s_{ip}(Y))^{1/p})\leq N/p$.
Plus g\'en\'eralement si $$(A_{m,n}(X,Y)=\suml{0\leq i\leq m/p}a_i(Y)X^i, \ B_{m,n}(X,Y)=\suml{1\leq j\leq m,\ (j,p)=1}b_j(Y)X^j),$$ on montre par r\'ecurrence sur $n$ que $i-w(a_i(Y))\leq N/p$ et que 
$j-w(b_j(Y))\leq N$. Ainsi $-\deg_Y(N_m(F^{[1]}(Y))=w(N_m(F^{[1]}(Y))=\suml{\tau\in G_n}w({F^{[1]}}^{\tau}(Y))\geq (1-N)p^{r(m)}$
et donc $\deg_Y(N_m(F^{[1]}(Y))\leq (N-1)p^{r(m)}=\deg \bar F'(Y)^{p^{r(m)}}$ puisque $(N,p)=1$. L'autre in\'egalit\'e est \'evidente. 

\Eproof
\subsection{R\'eduction stable}
\begin{Lem} .--- Le mod\`ele  $\cal C'$ de $C$, obtenu par normalisation  de $\Proj R[X_0,X_1]$, le mod\`ele de $\PP^1$ avec $X=\frac{X_1}{X_0}$, a une fibre sp\'eciale r\'eduite qui est un rev\^etement radiciel de 
$\PP^1_k$; les singularit\'es sont des cusps et se trouvent  au-dessus des z\'eros 
$\in \PP^1_k$ d'une forme  diff\'erentielle logarithmique $\omega$ r\'eguli\`ere  en dehors de $\bar \Br$. Le graphe d'intersection de la fibre sp\'eciale ${\cal C}_s$ du  mod\`ele stable $\cal C$ de la courbe point\'ee $(C,\Br)$  
est un arbre et le mod\`ele de $\PP^1$ obtenu par quotient de l'action de $\ZZ/p\ZZ$ sur $\cal C$ a une fibre sp\'eciale qui est un arbre de 
droites projectives attach\'ees  \`a $\PP^1_k$ en les points de $\bar \Br\subset \PP^1_k$.
\end{Lem}

\centerline {\input fig1.pstex_t}

\ms
{\it Preuve.} On consid\`ere le mod\`ele ${\cal C}'$ de $C$ obtenu par cl\^oture int\'egrale dans $K(C)$ du 
mod\`ele minimal de $\lv P^1$ qui d\'eploie le lieu de branchement; ici c'est le mod\`ele lisse 
de $\lv P^1$ correspondant \`a la coordonn\'ee $X$ de l'\'equation normalis\'ee $Z^p=F(X)$. La fibre 
sp\'eciale ${\cal C}'_s$ est hom\'eomorphe \`a $\lv P^1_k$; les singularit\'es sont des cusps.
Il reste \`a faire une \'etude locale: l'\'equation $Z^p=F(X)$ induit en r\'eduction une \'equation
$\bar Z^p=\bar F(X)$; il suit par le crit\`ere Jacobien que les cusps sont concentr\'es 
au-dessus des z\'eros de $\bar F'(X)$. Soit 
$x_1\in \Br$ et $\bar x_1$ sa sp\'ecialisation; une \'equation locale de ${\cal C}'$ dans la fibre 
formelle en $\bar x_1$ est ${Z'}^p=(X-x_1)\prodl{2\leq i\leq m}(X-x_i)^{u_1e_i}$ o\`u $u_1e_1=1 \mod p.$
Cette \'equation induit un mod\`ele lisse au-dessus de $\bar x_1$. On a donc montr\'e que les cusps ne 
sont pas au-dessus de  $\{\infty,\bar x_1,...,\bar x_m\}$ et se trouvent  au-dessus de z\'eros de  $\bar F'(X)$; plus 
simplement dit ils sont au-dessus des  z\'eros de $\omega =\frac{d\bar F}{\bar F}$  (cf. fig. 1). 
Notons $(\omega_0)$ (resp. $(\omega_\infty)$)
le diviseur des z\'eros  (resp. p\^oles) de $\omega$, alors $\deg (\omega_0)=\deg (\omega_\infty)-2=
m-1)$. Ainsi le 
mod\`ele stable est lisse au-dessus de $\lv P^1$ en dehors des z\'eros de $\omega$. Soit $x$ un tel 
z\'ero, apr\`es localisation en $x$ on se retrouve dans la situation locale d\'ecrite par Raynaud dans
([Ra1], d\'emonstration du th\'eor\`eme 1 p. 182), le lemme suit.  

Ce lemme montre que dans la situation o\`u le lieu de branchement est \'equidistant on aura un 
algorithme pour trouver la r\'eduction stable d\`es que l'on saura d\'etecter les composantes de genre non nul. Le th\'eor\`eme suivant donne un tel algorithme. 

\begin{The} Chaque composante de genre non nul du mod\`ele stable correspond \`a une valuation de Gauss sur un 
disque ferm\'e de $\lv P^1$ qui peut \^etre ainsi d\'efinie. Pour  $\bar d$ un z\'ero de 
$\bar F'(X)$ qui n'est pas un z\'ero de $\bar F(X)$, soit $m(\bar d)-1\leq m-1,$ son ordre  et soient 
$d_1,...,d_{\varphi(\bar d)}$
les z\'eros qui se sp\'ecialisent dans $\bar d $ de $N_m(F^{[1]}(Y))$, la $p$-d\'eriv\'ee, de $F(X)$ de niveau $m$, alors $\varphi(\bar d)=p^{r(m)}(m(\bar d)-1)$ (on a $\overline {N_m(F^{[1]}(Y))} =\bar F'(Y)^{p^{r(m)}}$). Pour $d_i$ avec $1\leq i\leq  \varphi(\bar d)$, consid\'erons $\tau_i\in G_n$ tel que ${F^{[1]}}^{\tau_i}(d_i)=0$, alors 
$$F(X)=E^{\tau_i}(X-d_i,d_i)^p+\suml{2\leq j\leq m,\, (j,p)=1}{F^{[j]}}^{\tau_i}(d_i)(X-d_i)^j \mod (p^{c_n}(X-d_i),(X-d_i)^{m+1}).$$
Soit $\rho_i \in K^{alg}$ tel que 
$$(*)\   v(\rho_i)=\maxl{ 2\leq j\leq m(\bar d),\, (j,p)=1} \frac{1}{j}v(\frac{\lambda^p}{{F^{[j]}}^{\tau_i}
(d_i)}),$$ 
alors la valuation de Gauss relative au disque $v(X-d_i)\geq v(\rho_i)$ induit dans le mod\`ele stable une composante de genre non nul; de plus toutes les composantes de genre non nul sont 
ainsi obtenues. Soit ${\cal P}^1$, le mod\`ele  semi-stable minimal de $(\PP^1,\Br)$ qui induit
ces valuations , alors le mod\`ele stable  $\cal C$ de $(C,G)$ est la cl\^oture int\'egrale de ${\cal P}^1$;  les bouts de l'arbre d'intersection repr\'esentent les composantes de genre non nul et l'origine
correspond \`a la composante  induite par le  disque $v(X)\geq 0$. Enfin les composantes de genre non nul sont des rev\^etements \'etales de la droite affine; en particulier la  jacobienne de $C$ a
potentiellement bonne r\'eduction supersinguli\`ere. 
\end{The}
{\it Preuve.} 
\begin{enumerate} [A.---] 
\item
Avant de passer \`a la preuve il faut remarquer que dans le cas o\`u $m<p$, le 
d\'eveloppement de Taylor classique suffit (2.1.2.i)); le
th\'eor\`eme est alors mot pour mot le th\'eor\`eme 4.1 de Lehr ([Le 2]). Passons \`a la preuve dans 
le cas g\'en\'eral; elle pr\'esente des difficult\'es nouvelles \`a cause de la complexit\'e des 
$p$-d\'eveloppements de Taylor lorsque $m\geq p$.  
\item 
Voyons que les valuations de Gauss relatives au disque $v(X-d_i)\geq v(\rho_i)$ induisent une 
composante de genre non nul dans le mod\`ele stable. 

Soit $m_i\geq 2$, $(m_i,p)=1$, le maximum des entiers $j,\ 2\leq j\leq m(\bar d)+1$ qui satisfont 
l'\'egalit\'e dans la relation $(*)$ du th\'eor\`eme. Notons que 
le $p$-d\'eveloppement de Taylor sp\'ecial induit modulo $(p)$ le d\'eveloppement de Taylor classique; ainsi 
par d\'efinition de $m(\bar d)$ on a $v(m(\bar d){F^{[m(\bar d)]}}^{\tau_i}(d_i))=0$ et donc 
$v(\rho_i)\geq \frac{1}{m(\bar d)}v(\lambda^p)\geq \frac{1}{m}v(\lambda^p)>0.$ On pose $X-d_i=\rho_iT$; l'\'equation du rev\^etement devient $$Z^p=F(X)=E^{\tau_i}(\rho_iT,d_i)^p+\suml{2\leq j\leq m,\, (j,p)=1}{F^{[j]}}^{\tau_i}(d_i){\rho_i}^j T^j\mod (p^{c_n}\rho_iT,(\rho_iT)^{m+1}).$$ 
Puisque $v(p^{c_n}\rho_i)\geq c_nv(p)+\frac{1}{m}v(\lambda^p)=(1-1/p^{n+1}+1/m)v(\lambda^p)>
v(\lambda^p)$ ({\bf c'est l\`a, la justification de l'approximation modulo $p^{c_n}$ dans les $p$-d\'eveloppements de Taylor d'ordre $m$}),
le changement $Z=\lambda^p W+E(\rho_iT,d_i)$ donne une \'equation enti\`ere qui en r\'eduction induit une \'equation s\'eparable 
$$\bar W^p-\bar W=\suml{2\leq j\leq m,\, (j,p)=1}\overline{{F^{[j]}}^{\tau_i}(d_i){\rho_i}^j/\lambda^p} \bar T^j=u_2\bar T+...+u_{m_i}\bar T^{m_i}$$ o\`u $u_{m_i}\neq 0$. On a ainsi exhib\'e une composante de ${\cal C}_s$ de genre $(m_i-1)(p-1)/2$ et d'invariant de Hasse-Witt nul. 
\Eproof

\item 
Passons au point le plus d\'elicat \`a savoir que l'on obtient bien ainsi toutes les valuations qui 
donnent un genre non nul dans le mod\`ele stable; il nous faut montrer que la somme des genres des 
composantes produites par l'algorithme comme au-dessus est $g(C)=(m-1)(p-1)/2$. Il n'est \`a priori pas garanti que ce soit le cas
car la condition d'annulation de $N_m({F^{[1]}}(Y))$ impos\'ee dans l'algorithme  est certainement 
restrictive; on pourrait trouver la m\^eme valuation de Gauss trop souvent et ainsi oublier des composantes. Nous allons voir qu'il n'en est 
rien.

Soit $N(d_i)$ le nombre (avec multiplicit\'e) de z\'eros $d$ de $N_m(F^{[1]}(Y))$ tels que $v(d-d_i)\geq v(\rho_i)$; nous allons montrer que $N(d_i)=p^{r(m)}(m_i-1)$. La strat\'egie est la suivante: 
puisque $N_m(F^{[1]}(Y))\in K[Y]$ alors pour $y\in R^{alg}$ avec
$v(y-d_i)<v(\rho_i)$ et proche de $v(\rho_i)$, il existe $C(d_i)\in \QQ$ tel que $v(N_m(F^{[1]}(y)))=
N(d_i)v(y)+C(d_i)$, une fonction lin\'eaire affine de $v(y)$. Nous allons utiliser deux $p$-d\'eveloppements de Taylor sp\'eciaux d'ordre $m$ de $F(X)$.
L'un centr\'e en $d_i$ et l'autre en $y$; ces deux d\'eveloppements donnent un mod\`ele pour l'extension valu\'ee au-dessus de la valuation de Gauss sur le disque ferm\'e 
$v(X-d_i)\geq v(y-d_i)$; sur chacun de ces mod\`eles on peut lire  la diff\'erente dans l'extension valu\'ee (d\'eg\'en\'erescence de $\mmu_p$ \`a $\aalpha_p$) et la comparaison des $\aalpha_p$-torseurs induits 
par chaque mod\`ele permet de conclure (cf. 1.1).  

Puisque pour $j\leq m(\bar d)$, $v({F^{[j]}}^{\tau_i}(d_i))+jv(\rho_i)\geq v({F^{[m_i]}}^{\tau_i}(d_i))+m_iv(\rho_i)=pv(\lambda)$, avec l'in\'egalit\'e sticte pour $m_i<j\leq m(\bar d)$; il suit que pour 
$\rho$ avec $0<v(\rho_i/\rho)<\epsilon$ et $\epsilon $ suffisament petit on conserve les m\^emes in\'egalit\'es si on remplace $\rho_i$ par $\rho$. Quitte \`a faire une extension des scalaires d'uniformisante encore not\'ee $\pi$, on peut supposer que 
$v({F^{[m_i]}}^{\tau_i}(d_i))+m_iv(\rho)=ptv(\pi)<pv(\lambda)$. Soit $Z=\pi^t W+E(X-d_i,d_i)$ et $X-d_i=\rho T$, puisque  $v(\lambda/\pi^t)>0$ 
alors $W^p=uT^{m_i} \mod \pi$, et si $w_{\frac{X-d_i}{\rho}}$ d\'esigne la valuation de Gauss associ\'ee au disque $v(X-d_i)\geq  v(\rho)$, la valuation de la diff\'erente dans l'extension valu\'ee 
correspondante est $v_\pi(p)-(p-1)t$ (cf. 1.1.1).

Pour $y$ avec $v(y-d_i)=v(\rho)$ et $\tau \in G_n$, on a un $p$-d\'eveloppement de Taylor sp\'ecial
 $$F(X=X-y+y)={E^{\tau}(X-y,y)}^p+\suml{(j,p)=1}{F^{[j]}}^{\tau}(y)(X-y)^j \mod (p^{c_n}(X-y),(X-y)^{m+1}).$$ Comme pr\'ec\'edemment on \'ecrit $\inf_j(v({F^{[j]}}^{\tau}(y))+jv(\rho))=pt'v(\pi)$.
Puisque $y$ est dans le disque $v(X-d_i)\geq v(\rho)$; la valuation de la diff\'erente dans l'extension valu\'ee par
$w_{\frac{X-d_i}{\rho}}$ est 
aussi $v_\pi(p)-(p-1)t'$  (cf. 1.1); ainsi $t=t'$. D'autre part si l'on pose $Z=\pi^tW'+E(X-y,y)^{\tau}$ 
et $X-y=\rho S$, il suit que $\quad {W'}^p\equiv c_0+c_1S+... \mod \pi$  qui est \`a comparer \`a $W^p\equiv uT^{m_i} \mod \pi$. On a $W-W'=\frac{E^{\tau}(X-y,y)-E(X-d_i,d_i)}{\pi^t}$ ainsi $\bar W^p-\bar {W'}^p\in (k[S])^p$ et donc  $uT^{m_i}\equiv u(S+(y-d_i)/\rho)^{m_i}\mod (\pi,S^p)$. Puisque $(m_i,p)=1$ et que $v((y-d_i)/\rho)=0$, il suit que $v(c_1)=0$ (comparer avec 1.1.1). Ainsi $v({F^{[1]}}^{\tau}(y))+v(\rho)=ptv(\pi)= m_iv(\rho)+v({F^{[m_i]}}^{\tau_i}(d_i))$; et donc  il existe  $C(d_i)$ avec $ v(N_m(F^{[1]}(y))=(m_i-1)p^{r(m)}v(\rho)+C(d_i)$; ainsi $N(d_i)=p^{r(m)}(m_i-1)$.
\Eproof

\item 
Soit $\bar d$, un z\'ero de $\omega=d\bar F/\bar F$ de multiplicit\'e $m(\bar d)-1$ 
et $d_i$, $d_j$ deux z\'eros de $N_m(F^{[1]}(Y))$ qui se sp\'ecialisent en  $\bar d$. 

Remarquons que la relation $v(d_i-d_j)\geq v(\rho_i)$ implique que $v(\rho_i)=v(\rho_j)$; en effet 
supposons que $v(d_i-d_j)\geq v(\rho_i)>v(\rho_j)$; on peut alors consid\'erer le $p$-d\'eveloppement 
de Taylor sp\'ecial d'ordre $m$ de  $F(X)$ centr\'e en $d_i$. Puisque $v(\rho_i)>v(\rho_j)>v(\lambda ^p)$, le changement de variable $X-d_i=\rho_j T$ induit un $\aalpha_p$ torseur en r\'eduction alors que par d\'efinition de $d_i$ et $\rho_i$ on a un $\ZZ/p\ZZ$-torseur! Contradiction.

Nous pouvons donc d\'efinir une partition index\'ee par $I(\bar d)$ des $\varphi(\bar d)$ z\'eros 
de $N_m(F^{[1]}(Y))$ qui se sp\'ecialisent en  $\bar d$ de la mani\`ere suivante: $d_i$ et $d_j$ seront \'equivalents si et seulement si $v(d_i-d_j)\geq v(\rho_i)=v(\rho_j)$. Notons que cette relation d'\'equivalence est la m\^eme que celle qui range les $d_i$ suivant la valuation de Gauss 
induite par le disque ferm\'e $v(X-d_i)\geq v(\rho_i)$. 
Chaque sous-ensemble de cette partition a un cardinal de la forme 
$p^{r(m)}(m_i-1)$ et induit une m\^eme composante dans le mod\`ele stable de genre $(m_i-1)(p-1)/2$. Puisque $\varphi(\bar d)=p^{r(m)}(m(\bar d)-1)=\suml{i\in I(\bar d)}p^{r(m)}(m_i-1)$, on obtient pour 
chaque $\bar d$ une contribution au genre dans le mod\`ele stable \'egale \`a $\suml{i\in I(\bar d)} (m_i-1)(p-1)/2=(m(\bar d)-1)(p-1)/2.$  Enfin puisque $\deg (\omega)_0=m-1$, $\suml{\bar d\in (\omega)_0}
(m(\bar d)-1)=m-1$; ainsi la somme des contributions au genre de ${\cal C}_s$ est $(m-1)(p-1)/2=g(C)$.

\Eproof 

\item
Puisque deux disques $v(X-d_i)\geq v(\rho_i)$ qui sont en inclusion sont confondus il suit qu'ils 
correspondent \`a des bouts dans l'arbre d'intersection (voir fig.2); ce dernier point est 
conforme \`a   [Ra 1], exemple (1) p. 186. Enfin la jacobienne de $C$  a potentiellement bonne r\'eduction supersinguli\`ere par   ([Ra 2]).

\Eproof
\end{enumerate}

\begin{Rem}.--- \begin{enumerate} [i)]{\rm 
\item  Soit $y\in R^{alg}$ tel que 
$v(N_m(F^{[1]}(y))\geq p^{r(m)}c_nv(p)$, il existe $\tau\in G_n$ avec $v(F^{[1]}(y))\geq c_nv(p)$, 
alors le $p$-d\'eveloppement de Taylor de $F(X+y)$ correspondant induit une composante de genre non nul en r\'eduction (m\^eme preuve que dans l'\'etape B); ainsi il existe un z\'ero $d$ de $N_m(F^{[1]}(Y))$ avec $v(d-y)\geq v(\rho_i)$. 
\item  La preuve montre que le nombre de composantes de genre non nul de la fibre sp\'eciale ${\cal C}_s$ du mod\`ele stable est major\'e  par $m-1$. Le polyn\^ome $N_m(F^{[1]}(Y)) $ est de 
degr\'e $\deg \bar F'(Y)^{p^{r(m)}}=(N-1) p^{r(m)}$ qui peut \^etre  tr\`es \'elev\'e puisque $p^{r(m)}$ est le 
cardinal du groupe de Galois de l'extension de $K(Y)$ engendr\'ee par les coefficients mis en jeu 
dans l'algorithme. Dans le cas o\`u $m<p$, on a vu en (2.1.2.i)) que $F'(Y)$ convient et puisque les 
seuls z\'eros \`a retenir sont ceux qui se sp\'ecialisent dans les z\'eros de 
$\omega =d\bar F/\bar F$, C. Lehr a introduit le polyn\^ome $N(Y)$ num\'erateur de la fraction irr\'eductible $F'(Y)/F(Y)$; ainsi 
$N(Y)$  divise $F'(Y)$, il  est de degr\'e $m-1$ et ses z\'eros donnent des centres pour les 
valuations correspondant aux composantes de genre non nul de la fibre sp\'eciale ${\cal C}_s$ du mod\`ele stable. Dans le cas g\'en\'eral il serait souhaitable d'exhiber un polyn\^ome de degr\'e minimal
sur $K$ ayant cette propri\'et\'e.} 
\end{enumerate}
\end{Rem}
\begin{Rem} .--- {\rm Il est possible de donner (cf. [He], [Sa 1,2,3]), des conditions combinatoires et diff\'erentielles qui sont des
conditions n\'ecessaires et suffisantes pour qu'une courbe stable sur $k$ soit de la forme ${\cal C}_s$
comme dans le th\'eor\`eme  (cf. fig. 2) }.
\end{Rem}
\centerline{\input fig2.pstex_t}
\subsection
{Bonne r\'eduction}
 \indent{\ \ \ }
Il suit facilement du th\'eor\`eme 3.2.2 un crit\`ere de potentielle bonne r\'eduction. 

\begin{The} .--- La courbe $C$ a potentiellement bonne r\'eduction si et seulement si le 
mod\`ele stable $\cal C$ de $(C,G)$ a 2 composantes \`a savoir  une composante de genre $g(C)$ et une droite projective sur laquelle se sp\'ecialise le lieu de branchement en $m+1$ points distincts lisses. C'est le cas si et seulement si 
 les z\'eros $y_i$ de $N_m(F^{[1]}(Y))$, la $p$-d\'eriv\'ee de niveau $m$ de $F(X)$, qui ne se sp\'ecialisent pas dans les z\'eros de $F(X)$ sont tels que $v(y_i-y_j)\geq \frac{1}{m}v(\lambda^p).$  
\end{The}
{\it Preuve.} En effet il suffit de reprendre le d\'ebut de la preuve du th\'eor\`eme 3.2.2. On a 
potentiellement bonne r\'eduction si et seulement un z\'ero de $N_m(F^{[1]}(Y))$ qui ne se sp\'ecialisent pas dans un z\'ero de $F(X)$ donne naissance \`a une composante de genre $(m-1)(p-1)/2$; i.e. de 
conducteur $m_i=m$. Ceci ne se produit que si $v(\rho_i)=\frac{1}{m}v(\lambda^p)$ et si les racines 
de $N_m(F^{[1]}(Y))$ donnent naissance \`a la m\^eme valuation.
\Eproof

Ce crit\`ere n\'ecessite  le calcul de $N_m(F^{[1]}(Y))$ et une localisation non triviale de ses
racines. Dans un cas particulier
C. Lehr a donn\'e un crit\`ere simple, nous allons rappeler ce crit\`ere
et voir qu'il d\'ecoule facilement de notre algorithme.

\begin{The}{\bf(\cite{Le2}, th\'eor\`eme 3.1)}.---
 Soit $F(X)$ comme au-dessus. Ecrivons $\frac{F'(X)}{F(X)}=\frac{N(X)}{D(X)}$ o\`u  $N(X)$ est unitaire  et  $(N(X),D(X))=1$. 
\begin{enumerate}[i)] 
\item  Si  $C$ a potentiellement bonne r\'eduction, alors  $(m,p)=1$ et il existe $d\in \tilde R$ tel que  $N(X)=(X-d)^{m-1} \mod (\lambda ^{p/m})$.
\sms
\item  Si $(m-1,p)=1$, si $C$ a  bonne r\'eduction, alors la  congruence ci-dessus d\'etermine
la classe de $d \mod \lambda^p$; precis\'ement  $-(m-1)d=-\suml{i} x_i+\frac{1}{\suml {i} e_i}\suml {i} e_ix_i\in R$ et la valuation de Gauss relative \`a  $\frac{X-d}{\lambda^{p/m}} $ induit le mod\`ele lisse
sur $R[f(d)^{1/p}]$. 
\end{enumerate}  
\end{The}
{\it Preuve.} Nous remontrons la congruence en utilisant notre algorithme. Le torseur en r\'eduction 
relativement \`a $v_X$, la valuation de Gauss associ\'ee \`a $X$, ne pr\'esente qu'un cusp (cf. lemme 3.2.1); ainsi $\bar F'(X)/\bar F(X)$ n'a qu'une racine qui est donc d'ordre $m-1$; en particulier $(m,p)=1$.
Soit $y$ une racine de $N_m(F^{[1]}(Y))$, alors $v(F^{[m]}(y))=0$ et donc  $\rho=\lambda^{p/m}$. 
Le $p$-d\'eveloppement de Taylor sp\'ecial d'ordre $m$ donne 
$F(X)=E(X-y,y)^p+F^{[1]}(y)(X-y)+...+F^{[m]}(y)(X-y)^m \mod (p^{c_n}(X-y),(X-y)^{m+1})$ et 
$v(\rho)\geq \frac{1}{j}(v(\lambda^p)-v(F^{[j]}(y)))$ pour $(j,p)=1$ et $j\leq m$, avec \'egalit\'e pour 
$j=m$; ainsi 
$v(F^{[j]}(y))\geq (1-j/m)v(\lambda^p)\geq v(\lambda ^{p/m})$ pour $j<m$ et premier \`a $p$.
On a $F(X)=\prodl{1\leq i\leq  m} (X-x_i)^{e_i} $ et donc avec  $N=\sum e_i$
$$N(X)= \frac{1}{N} \frac{F'(X)}{F(X)}\prodl{1\leq i\leq m}(X-x_i);$$ 
ainsi 
$$N(X+y)=$$
$$=\frac{1}{N}[pE(X,y)^{p-1}E'(X,y)+\suml{(j,p)=1}F^{[j]}(y)jX^{j-1}
+\lambda^{p-1/p^n}C(X)+ X^mD(X)]\prodl{1\leq i\leq m}(y-x_i+X)^{-e_i+1}=$$ 
$$=X^{m-1}\prodl{1\leq i\leq m}(1+X/(y-x_i))^{-e_i+1}\mod (\lambda ^{p/m},X^m)R[[X]]$$ qui donne la congruence du th\'eor\`eme puisque $N(X)$ est 
un polyn\^ome de degre $m-1$. 

Le reste de la preuve est imm\'ediat. 

\Eproof
\begin{Rem}.--- \begin{enumerate}[i)]{\rm 
\item  Si on a un centre $y$ et un rayon candidats il est facile de tester la bonne r\'eduction; cf. la proposition 2.5 de [Le 2].
\item  La condition $(m-1,p)=1$ est tr\`es restrictive en particulier pour les petits
$p$ et si $p=2$; le crit\`ere de Lehr est vide.
\item  On peut se demander si il n'y a  pas d'autres cas o\`u apparait un bon candidat pour
le centre par exemple on aimerait  preciser un centre pour chaque composante
de la r\'eduction stable de conducteur $m_i$ avec $(m_i-1,p)=1$? Supposons que $F'$ a 
un seul z\'ero en $y$ d'ordre $m_y-1$ avec $(m_y-1,p)=1$; on veut alors tester si il n'y a qu'une composante dans la r\'eduction stable au dessus du cusp correspondant a $y$. Dans ce cas il n'y a qu'un rayon
candidat qui est $v(\lambda ^{p/m_y})$; la preuve pr\'ec\'edente montre alors si l'on appelle $\tilde y$ un 
z\'ero de $N_m(F^{[1]}(Y)) $ qui se sp\'ecialise en $y$ on a la congruence 
$N(X)=(X-\tilde y)^{m_y-1}A(X) \mod (\lambda ^{p/m_y})$ pour $A(X)=X^{m-m_y}+\alpha_1X^{m-m_y-1}+...\in R[X]$. Il suit que  $-(m_y-1)y+\alpha_1$ est connu modulo $\lambda ^{p/m_y}$; mais 
c'est insuffisant pour conclure...}
\end{enumerate}
\end{Rem}
\section {Exemples}
\indent{\ \ \ }
Nous reprenons les notations de la partie 3.

Dans cette partie nous illustrons les th\'eor\`emes pr\'ec\'edents. Comme dans ([Ra1]), la question 
de la situation g\'en\'erale se pose. Nous passerons ensuite \`a des exemples plus sp\'ecifiques.
\subsection{Le cas de ramification g\'en\'erale}
\indent{\ \ \ }
Nous imposons des conditions g\'en\'erales au lieu de branchement (toujours sous l'hypoth\`ese 
d'\'equidistance); en fait la condition est que la sp\'ecialisation du lieu de branchement est en position 
g\'en\'erale; la proposition suivante est de ce point de vue un analogue facile de ([Ra1], cor. 4 p. 194). 

\begin{Pro} 
.---  Soit $Z^p=F(X)=\prodl{1\leq i\leq m}(X-x_i)^{e_i}$ avec $(e_i,p)=1$ et 
$(\deg F,p)=1$. On suppose que $v(x_i-x_j)=v(x_i)=0$ pour tout $i\neq j$. Alors il existe un polyn\^ome
$\Delta (X_1,...,X_m)\in k[X_1,...,X_m]-\{0\}$ tel que la condition $\Delta (\bar x_1,...,\bar x_m)\neq 0$ implique que la forme diff\'erentielle $\omega=d\bar F(X)/\bar F(X)$ n'a que des z\'eros simples 
$z_1,...,z_{m-1}$. Ainsi la
fibre sp\'eciale du mod\`ele stable ${\cal C}$ de $(C,G)$ qui d\'eploie le lieu de branchement est 
constitu\'ee d'une part d'une  droite projective sur $k$ contenant  les $m+1$ sp\'ecialisations  du lieu de branchement et les $m-1$ z\'eros de $\omega $; et d'autre part de $m-1$ composantes de genre $(p-1)/2$ si $p>2$ et de genre $1$ si $p=2$ intersectant la droite projective dans les z\'eros de $\omega $.
\end{Pro}
 
{\it Preuve.} Afin de montrer l'existence de $\Delta $ nous montrons 2 lemmes suivant que $p>2$ 
ou
que $p=2$. 
\begin{Lem}.--- 
Soit $p>2$, $m\geq 3$ et $f(X):=\prodl{1\leq i\leq m}(X-x_i)^{e_i}\in k[X]$
o\`u $(e_i,p)=1$ et $(N=\suml{i}e_i,p)=1$. Soit  $S_m(x_i,e_i)(X):=
\suml{1\leq i\leq m}e_i\prodl{1\leq j\leq m,j\neq i}(X-x_j)$ alors son discriminant $\disc (S_m(x_i,e_i)(X))\in \lv F_p[x_1,x_2,...,x_m]$  n'est pas identiquement nul.
\end{Lem}
{\it Preuve.} Nous faisons une preuve par r\'ecurrence sur $m$.  Si $m=3$, $x_1=0,x_2=1$, alors 
 $$S_3(x_i,e_i)(X)=e_3X(X-1)+e_2X(X-x_3)+e_1(X-1)(X-x_3)=$$
$$=(e_1+e_2+e_3)X^2+(-e_3-e_2x_3-e_1-e_1x_3)X+e_1x_3$$ son discriminant est $$(-e_3-e_2x_3-e_1-e_1x_3)^2-4(e_1+e_2+e_3)e_1x_3= $$
$$=(e_1+e_2)^2x_3^2+
(2(e_1+e_3)(e_1+e_2)-4e_1(e_1+e_2+e_3)x_3+(e_1+e_3)^2$$ qui n'est pas identiquement nul.

Pour $m>3$ nous utilisons un argument de sp\'ecialisation: soit $f(X)$ comme dans le lemme, on 
peut supposer que $e_{m-1}+e_m\neq 0 \mod p$, alors $$S_m(x_i,e_i)(x_m=x_{m-1})(X)=(X-x_{m-1})S_{m-1}((x_1,e_1),...,(x_{m-2},e_{m-2}),(x_{m-1},e_{m-1}+e_m))(X);$$ de plus $$S_{m-1}((x_1,e_1),...,
(x_{m-2},e_{m-2}),(x_{m-1},e_{m-1}+e_m))(X=x_{m-1})
=(e_{m-1}+e_m)\prodl{1\leq i\leq m-2}(x_{m-1}
-x_i)$$ et par r\'ecurrence, pour $x_1,...,x_{m-1}$ en position g\'en\'erale $S_m(x_i,e_i)(x_m=x_{m-1})(X)$ est de degr\'e $m-1$ et n'a que des racines simples; ainsi $\disc (S_m(x_i,e_i)(X))$, n'est pas
identiquement nul.
\Eproof
\begin{Lem}.--- 
Si $p=2$, soit $f(X):=\prodl{1\leq i\leq 2m+1}(X-x_i)\in k[X]$, avec  
$x_1,...,x_{2m+1}\in k^{2m+1}$, alors $f'(X)=S(X)^2$ et $\disc f(X)\disc S(X)\in \lv F_2[x_1,x_2,...,x_m]$ n'est pas identiquement nul. 

\end{Lem}
{\it Preuve.} Soit $S(X)=\prodl{1\leq i\leq m}(X-y_i)\in k[X]$ avec $y_i\neq y_j$ si $i<j$. Soit 
$f(X)=1+S(X)^2+XS(X)^2$ alors  $f'(X)=S(X)^2$. Ainsi $\disc f(X)\disc S(X)\neq 0$. 
\Eproof

Nous passons \`a la preuve de la proposition.

Si $p>2$, le lemme 4.1.2  montre qu'il existe $\Delta (X_1,...,X_m)\in k[X_1,...,X_m]-\{0\}$ tel que la condition $\Delta (\bar x_1,...,\bar x_m)\neq 0$ implique que la forme diff\'erentielle $\omega=d\bar F(X)/\bar F(X)$ dont le num\'erateur vaut $S_m(x_i,e_i)(X)$ n'a que des z\'eros simples 
$z_1,...,z_{m-1}$; on se retrouve localement avec un conducteur $<p$. Conform\'ement \`a [Le 2], les 
z\'eros $y$ de la d\'eriv\'ee $F'$ tels que $F(y)$ est une unit\'e  fournissent des centres des disques qui induisent un genre non nul dans
le mod\`ele stable. Alors pour un tel z\'ero on a
$F(X+y)=F(y)+\frac {F"(y)}{2} X^2+...$, il suit du lemme 4.1.2 que  $v(F"(y))=0$. 
Soit $X=\lambda^{p/2}T$ alors  $F(y+\lambda^{p/2}T)-F(y)\in \lambda^{p}R^{alg}[T]$ induit un 
rev\^etement \'etale de la droite affine de conducteur $m_i=2<p$ comme annonc\'e (cf. fig. 3)
\vskip 40pt
\centerline{\input fig3.pstex_t}
\vskip 15pt
Si $p=2$, le lemme 4.1.3. montre qu'il existe $\Delta (X_1,...,X_m)\in k[X_1,...,X_m]-\{0\}$ tel que la condition $\Delta (\bar x_1,...,\bar x_m)\neq 0$ implique que la forme diff\'erentielle $\omega=d\bar F(X)/\bar F(X)$ n'a que des z\'eros 
$z_1,...,z_{(m-1)/2}$ de multiplicit\'e $2$; contrairement au cas pr\'ec\'edent on se retrouve localement avec un conducteur $>p=2$. Nous appliquons l'algorithme du th\'eor\`eme 3.2.2. Consid\'erons $y\in  \tilde R$ tel que 
$\bar y$ soit racine de $\bar F '(X)$; nous la choisirons plus pr\'ecis\'ement plus tard. Ecrivons 
le $2$-d\'eveloppement de Taylor sp\'ecial de $F$ de niveau $3$, on a
$F(X+Y)=s_0(Y)+s_1(Y)X+s_2(Y)X^2+ s_3(Y)X^3...=(s_0(Y)^{1/2}+s_2(Y)^{1/2}X)^2
+(s_1(Y)-2s_0(Y)^{1/2}s_2(Y)^{1/2})X+ s_3(Y)X^3...$; la norme du coefficient de $X$ donne
 $N_3(F^{[1]}(Y))=s_1(Y)^2-4s_0(Y)s_2(Y)={F'}^2-2FF"$. Soit $y$ une racine de $N_3(F^{[1]}(Y))$ qui 
se r\'eduit sur une racine $\omega=d\bar F/\bar F$; alors $s_0(y)=f(y)$ et $s_3(y)$ sont des unit\'es
et pour un choix convenable de la racine $y$ on obtient: $f(X+y)= (s_0(y)^{1/2}+s_2(y)^{1/2}X)^2 +s_3(y)X^3...$. Soit $X=\lambda ^{2/3}T$ et $Z=\lambda W+(s_0(y)^{1/2}+s_2(y)^{1/2}X)$, alors en r\'eduction 
on obtient l'\'equation $W^2-W=\bar s_3(y)T^3 $. Pour conclure il n'est nul besoin d'invoquer le 
th\'eor\`eme 3.2.2 puisque chaque z\'ero de  $\omega=d\bar F(X)/\bar F(X)$ donne naissance 
\`a une (seule) composante; il est alors imm\'ediat de v\'erifier que l'on a le bon genre \`a la
fibre sp\'eciale. (cf. fig. 4)
\centerline{\input fig4.pstex_t}

\subsection{Les courbes hyperelliptiques en $p=2$}
\indent{\ \ \ }
On a donc $F(X)=s_0+s_1X+....+X^m$ avec $(m,2)=1$ et $\disc  \bar F\neq 0$. Nous allons examiner le cas o\`u $m$ est petit. 
\subsubsection { L'algorithme de Coleman revisit\'e}
\indent{\ \ \ }
Dans le cas de petites valeurs de $m$ on peut tenter d'obtenir un $2$-d\'eveloppement de Taylor exact de 
$F(X)$ de niveau $m$  (cf. 2.4.1). Pour $m\leq 7$, nous \'ecrivons les \'equations de la vari\'et\'e de dimension $0$ correspondante et calculons la $2$-d\'eriv\'ee qui lui est attach\'ee en calculant
comme dans le cas des d\'eveloppements sp\'eciaux. Enfin nous suivons notre algorithme pour calculer 
la $2$-d\'eriv\'ee de niveau $m$; les deux m\'ethodes donnent des $2$-d\'eriv\'ees
comparables. Pour $m>9$, la m\'ethode de Coleman donne naissance \`a une vari\'et\'e dont nous ne connaissons pas la non vacuit\'e et les calculs via notre algorithme deviennent vite compliqu\'es.

\begin{enumerate}[A.]
\item {\bf $m=3$ i.e. $g=1$.}

On a donc $F(X+Y)=s_0(Y)+s_1(Y)X+s_2(Y)X^2+X^3=(a_0+a_1X)^2+ b_1X+X^3$  que l'on r\'esout en 
$a_0,a_1,b_1$; on a $a_0^2=s_0(Y)$, $a_1^2=s_2(Y)$, $b_1=s_1(Y)\pm 2s_0(Y)^{1/2}s_2(Y)^{1/2}$, ainsi 
la norme dans un sens \'evident de $b_1$ est $N_3(b_1)=s_1(Y)^2-4s_0(Y)s_2(Y)= {F'}^2-2FF"=
-3 Y^4-4 s_2 Y^3-6 s_1 Y^2-12 s_0 Y+s_1^2-4 s_0 s_2=
F'(Y)^2 \mod 2$. 

Appliquons l'algorithme: $F(X+Y)=(s_0(Y)^{1/2}+s_2(Y)^{1/2}X)^2+(s_1(Y)-2s_0(Y)^{1/2}s_2(Y)^{1/2})X+X^3$;
c'est un d\'eveloppement exact; ainsi la $2$-d\'eriv\'ee de niveau $m$ est aussi $N_3(F^{[1]}(Y))={F'}^2-2FF"$. 

Soit $y$ une racine de $N_3(F^{[1]}(Y))$ alors la valuation de
Gauss relative au disque $v(X-y)\geq \frac{2}{3}v(2)$ induit une courbe elliptique en r\'eduction. On a
potentiellement bonne r\'eduction.

\item
{\bf $m=5$ i.e. $g=2$. } 
On a cette fois 
$(*)\quad F(X+Y)=s_0(Y)+s_1(Y)X+s_2(Y)X^2+s_3(Y)X^3+s_4(Y)X^4+X^5=(a_0+a_1X+a_2X^2)^2+ b_1X+b_3X^3+X^5$
  
que l'on r\'esout en $a_0,a_1,a_2,b_1,b_3$. On a

$$s_0(Y)=a_0^2$$
$$s_1(Y)=2a_0a_1+b_1$$
$$s_2(Y)=2a_0a_2+a_1^2$$
$$s_3(Y)=2a_1a_2+b_3$$
$$s_4(Y)=a_2^2$$

On utilise le syt\`eme de calcul Maple:

$b_1:=s1-2*RootOf(t^2-s0,t)*RootOf(t^2-(s2-2*RootOf(u^2-s0,u)*RootOf(u^2-s4,u)),t);$
la commande $evala(Norm(b_1));$ donne 

$N(b1)=(-64 s_4 s_0^3  + 16 s_0^2  s_2^2  - 8 s_0 s_2 s_1^2  + s_1^4 )^2=((s_1(Y)^2-4s_0(Y)s_2(Y))^2-64s_0(Y)^3s_4(Y))^2,$ (le carr\'e vient du passage \`a une cl\^oture galoisienne).

Appliquons l'algorithme, on a cette fois-ci 

$F(X+Y)=(s_0(Y)^{1/2}+s_2(Y)^{1/2}X+s_4(Y)^{1/2}X^2)^2+(s_1(Y)-2s_0(Y)^{1/2}s_2(Y)^{1/2})X \linebreak -2s_0(Y)^{1/2}s_4(Y)^{1/2}X^2 -2s_2(Y)^{1/2}s_4(Y)^{1/2}X^3+X^5;$

 ce n'est pas un 2-d\'eveloppement de niveau 5; cependant l'\'etape suivante modifie le coefficient de $X$ par
 $-2s_0(Y)^{1/2}(-2s_0(Y)^{1/2}s_4(Y)^{1/2})^{1/2}$. Ainsi
la $2$-d\'eriv\'ee de niveau $5$ est $N_5(F^{[1]}(Y)) =Norm((s_1(Y)-2s_0(Y)^{1/2}s_2(Y)^{1/2}-2s_0(Y)^{1/2}(-2s_0(Y)^{1/2}s_4(Y)^{1/2})^{1/2})$ que l'on calcule avec Maple 

 $N_5(F^{[1]}(Y)):=
(-3072s_4s_0^4s_2s_1^2+4096s_4^2s_0^6-2048s_2^2s_4s_0^5-128s_4s_0^3s_1^4+256s_0^4s_2^4-
256s_0^3s_2^3s_1^2+96s_1^4s_2^2s_0^2-16s_1^6s_2s_0+s_1^8)^2$

Cette norme est tr\`es proche du carr\'e de la pr\'ec\'edente. Ceci illustre le fait que parmi les 
$2$-d\'eveloppements de niveau donn\'e les sp\'eciaux ne sont pas n\'ecessairement les meilleurs.
Dans le cas pr\'esent la m\'ethode de Coleman est plus avantageuse. 

Le th\'eor\`eme 3.2.2 s'applique au $2$-d\'eveloppement (*); ainsi   
soit $y\in R^{alg}$ une racine de $N(b_1)$ qui r\'esiduellement est racine de 
$\omega =d\bar F/\bar F$  ou ce qui revient au m\^eme   de $\bar F'(Y)$ (puisque $p=2$). 
En particulier $s_0(y)$ est une unit\'e; on peut supposer que $F^{[1]}(y)=0$.  Suivant le 
th\'eor\`eme II. 2.2 on calcule 
$$(**)\   v(\rho)=\max(\frac{1}{3}v(\frac{2^2}{b_3(y)}),\frac{1}{5}v(2^2)).$$ 
On distingue alors 3 cas de figure:

\begin{enumerate}[1.---]
\item {\bf $\bar F'(Y)$ a un seul z\'ero qui est d'ordre $4$ et $C$ a potentiellement bonne r\'eduction.}
Ce cas se produit si et seulement si  l'\'egalit\'e dans la formule (**) est r\'ealis\'ee pour $j=5$
i.e. $v(\rho)=\frac{1}{5}v(2^2)$.
Par exemple c'est le cas pour $F(X)=1+X^5$ (\'evident) ou pour $F(X)=1+2X+X^5$ et dans ce n'est pas \'evident (on peut v\'erifier que c'est bien le cas en calculant les invariants d'Igusa avec l'algorithme de Liu, cf. [Li 2]). 
On calcule $N(b_1(y))=-320*y+1600*y^6-320*y^3-1600*y^7-2000*y^{11}+16+160*y^4+600*y^8+1000*y^{12}+625*y^{16}$, un algorithme utilisant le syst\`eme de calcul PARI montre que $\inf v(y_i-y_j)=(2/5)v(2)$. On a donc potentielle bonne r\'eduction. 

\item {\bf $\bar F'(Y)$ a un seul z\'ero qui est d'ordre $4$ et $C$ n'a pas potentiellement bonne r\'eduction.} 
L'existence d'une telle configuration (cf. fig. 5) est garantie par le recollement formel d\`es que les 
conditions diff\'erentielles \'evoqu\'ees en 3.2.4 sont satisfaites: dans la figure 5 consid\'erons
la composante $\PP^1_k$ horizontale; le torseur d\'eg\'en\`ere en $\aalpha_2$ sur cette composante et
les 2 courbes elliptiques $E$ et $E'$ intersectent cette droite
en deux points $p$ et $p'$ qui correspondent aux z\'eros de la forme diff\'erentielle additive $\omega_1=t^2(1+t)^2dt$ dont le p\^ole est le point d'intersection avec la composante initiale. 

Ce cas de figure se produit par exemple pour $F(X)=1+2^{1/2}X^3+X^5$ puisque le changement $X=2^{1/2}T$
 induit une composante de genre 1. D'autre part le changement $X=2^{1/4}T$ donne $Y^2= 1+2^{5/4}(T^3+T^5)$ qui induit un $\aalpha_2$ torseur $Y^2=t^3+t^5:=f$ et $df=t^2(1+t)^2dt$. La d\'etermination de l'autre composante de genre 1 n\'ecessite l'algorithme: on calcule $N(b1(y))=-95*Y^{16}-300*2^{1/2}*Y^{14}-772*Y^{12}+240*Y^{11}-376*2^{1/2}*Y^10+200*2^{1/2}*Y^9+36*Y^8-384*Y^7+640*Y^6+144*2^{1/2}*Y^5+288*Y^2-320*Y=Y^{16}+4*2^{1/2}*Y^{14}+4*Y^{12}+4*Y^8 \mod 2^3$; ainsi si $y$ est une racine de valuation $(1/4)v(2)$, puisque 
$b_3(y)=2^{1/2} \mod 2$ il suit que le changement $X=y+\rho T$ avec $v(\rho)=(1/2)v(2)$ induit la 
deuxi\`eme composante de genre 1. 

Consid\'erons maintenant le cas de  $F(X)=1+2X+2^{1/2}X^3+X^5$ . Contrairement au cas pr\'ec\'edent il n'y a pas de z\'ero \'evident pour 
$N(b_1(y))=(3084*y^8+1600*y^6-4480*y^7-2000*y^{11}+3700*y^{12}+625*y^{16}-96*2^{1/2}*y+96*2^{1/2}*y^2+16-1872*2^{1/2}*y^5+960*2^{1/2}*y^4-3000*2^{1/2}*y^9+1152*2^{1/2}*y^6+2880*2^{1/2}*y^{10}+1500*2^{1/2}*y^{14}-320*y+592*y^4-896*y^3+288*y^2)^2$ cependant l'analyse du polygone de Newton montre 
que les racines du  polyn\^ome sont de valuation $\frac{1}{4}v(2)$ qui se r\'epartissent en 2 classes
modulo $2^{1/4}$; ainsi on peut trouver 2 racines $y_1$ et $y_2$ avec $v(y_1-y_2)=v(y_1)=\frac{1}{4}v(2)<\frac{2}{5}v(2)$ ce qui contredit le crit\`ere de bonne r\'eduction (cf. Th. 3.3.1).   
\vskip 10pt

\centerline{\input fig5.pstex_t}

\item {\bf $\bar F'(Y)$ a deux z\'eros qui sont d'ordre $2$.}
On trouve 2 composantes de genre $1$ mais contrairement au cas pr\'ec\'edent c'est la configuration
g\'en\'erique (cf. fig. 4).
\end{enumerate}
\item
{\bf $m=7$ i.e. $g=3$.} 
On distingue alors 6  cas de figure qui sont num\'erot\'es en fonction de la g\'eom\'etrie diff\'erentielle sous-jacente en respectant l'application de sp\'ecialisation dans un certain espace de modules. Par exemple $\Delta_0$ correspond au cas o\`u la forme diff\'erentielle logarithmique $\omega$ a 3 zeros
de multiplicit\'e 2 (c'est le cas g\'en\'eral vu en 4.1). Le cas $\Delta_{010}$ correspond au cas o\`u la forme diff\'erentielle logarithmique $\omega$ a 2 zeros dont l'un est de multiplicit\'e 2 et l'autre de multiplicit\'e 4. Ce dernier cas
se sp\'ecialise en $\Delta_{010}$ et correspond au cas o\`u la forme diff\'erentielle additive sur la
droite projective sur laquelle s'appuient les 2 courbes elliptiques a 2 poles d'ordre 2; ce cas lui-m\^eme se sp\'ecialisant dans $\Delta_{011}$. 

On peut verifier que les conditions diff\'erentielles sont r\'ealisables et donc via le recollement 
formel on peut construire des exemples correspondant \`a chacune des 6 situations (cf. fig. 6).
Nous allons donner un exemple num\'erique sous la forme d'une \'equation globale $Y^2=F(X)$ dans chacun des cas. Comme dans le cas du genre 1 on se simplifie la vie en choisissant un exemple pour lequel une
valuation qui induit une composante de genre non nulle correspond \`a un disque ferm\'e centr\'e en $X=0$; on a ainsi $Y=0$ qui est un z\'ero \'evident de la $7$-d\'eriv\'ee.

On applique la m\'ethode de Coleman; cette fois 

$F(X+Y)=s_0(Y)+s_1(Y)X+s_2(Y)X^2+s_3(Y)X^3+s_4(Y)X^4+s_5(Y)X^5+s_6(Y)X^6+X^7=(a_0+a_1X+a_2X^2+a_3X^3)^2+ b_1X+b_3X^3+b_5X^5+X^7$

  que l'on r\'esout en $a_0,a_1,a_2,a_3,b_1,b_3,b_5$. On a

$s_0(Y)=a_0^2,\ s_1(Y)=2a_0a_1+b_1,\ s_2(Y)=2a_0a_2+a_1^2,\ s_3(Y)=2a_0a_3+2a_1a_2+b_3,\ s_4(Y)=2a_1a_3+a_2^2,\ s_5(Y)=2a_2a_3+b_5,\ s_6(Y)=a_3^2.$

En utilisant le syt\`eme de calcul Maple on obtient 

$N(b_1):=(-4096*s6*s0^5*s1^2+256*s0^4*s2^4-2048*s0^5*s2^2*s4-256*s0^3*s2^3*s1^2+96*s0^2*s2^2*s1^4+4096*s0^6*s4^2+1024*s0^4*s4*s1^2*s2-128*s0^3*s4*s1^4-16*s0*s1^6*s2+s1^8)^2$

Nous utilisons l'algorithme pour exhiber $N_7(F^{[1]}(Y))$, la $2$-d\'eriv\'ee de niveau 7. 
Suivant le m\^eme principe que pour $m=5$, apr\`es 2 \'etapes on obtient comme coefficient de $X$;  
$s_1(Y)-2*(-2)^{1/2}*((s_0(Y))^{1/2}*(s_4(Y))^{1/2})^{1/2}*(s_0(Y))^{1/2}-2*(s_0(Y))^{1/2}*(s_2(Y))^{1/2}$; l'\'etape suivante le modifie par un terme de taille $2^{1+1/2+1/4}$ et puisque $1+1/2+1/4+2/7>2$ on peut
le n\'egliger pour exprimer la $2$-d\'eriv\'e de niveau $7$; ainsi 

$N_7(F^{[1]}(Y))=(4096*s0^6*s4^2-128*s0^3*s4*s1^4-3072*s0^4*s4*s1^2*s2-2048*s0^5*s4*s2^2+s1^8-16*s1^6*s2*s0+96*s1^4*s2^2*s0^2-256*s1^2*s2^3*s0^3+256*s2^4*s0^4)^8$

les radicaux de $N(b_1)$ et $N_7(F^{[1]}(Y))$ sont de degr\'e 48 en $Y$ et  sont \'egaux modulo 
$2^{12}$. 

\quad

\quad 

\centerline{\input fig6.pstex_t}

\begin{enumerate}[1.---]
\item {\bf $\Delta_{0}$}

On prend $F(X)=1+X^3+X^5+X^7$; alors $\bar F'(X)=(X+X^2+X^3)^2$ et $X+X^2+X^3$ a trois racines simples.
C'est la situation g\'en\'erale. On calcule $N:=N_7(F^{[1]}(Y))$ d'o\`u les congruences 
$\rem(N,y^9,y) \mod 2^11=1024*y^6+256*y^7+256*y^4$ et 
$N mod 2=
y^32+y^16+y^48$; les racines sont regroup\'ees dans 3 classes modulo $2$ et si $y$ est une racine on a
$F^{[3]}(y)=1+y^4 \mod 2$ ainsi $v(F^{[3]}(y))=0$ et le changement $X=y+2^{2/3}T$ induit une courbe de genre 1.

\item {\bf $\Delta_{010}$}

On part du polyn\^ome  $\bar F(X)=1+X^5+X^7$ dont la diff\'erenrielle $X^4(1+X)^2dX$ a 2 z\'eros d'ordre respectifs $4$ et $2$; on va montrer que $F(X)=1+2^{1/2}X^3+X^5+X^7$ convient. Remarquons que le changement $X=2^{1/4}T$ induit une equation $Y^2=1+2^{5/4}(T^3+T^5)+2^{7/4}T^7$ qui induit un $\aalpha_2$-torseur $Y^2=T^3+T^5:=f$ et $df=(T+T^2)^2dT $ a 2 racines doubles; on trouve (certainement) ainsi la composante qui croise les 2 courbes elliptiques. Pour
confirmer cela, on calcule la $7$-d\'eriv\'ee $N:=N_7(F^{[1]}(Y))$ et on analyse son polygone de Newton, quelques congruences permettent de conclure: $ N=Y^{32}+Y^{48} \mod 2$, $\rem (N,Y^{33},Y) \mod 2=Y^{32}$
$\rem(N,Y^3,Y) \mod 2^{13}=4096*Y^2$, $\rem(N,Y^{17},Y) \mod 2^5= 16*Y^{16}$; ainsi il y a 16 racines de valuation $\geq (4/7)v(2)$; 16 autres de valuation $(1/4)v(2)$ et enfin les 16 autres de valuation nulle. Le 
premier lot correspond  \`a la valuation sur le disque $v(X-y)\geq 1/2$ qui induit une courbe elliptique (notez que $(4/7)v(2)>(1/2)v(2)$).

\item {\bf $\Delta_{011}$}

On prend $F(X)=1+X^5+X^7$; alors $(\bar F(X))'=X^4(1+X)^2$. Le changement $X=2^{2/5}T$ induit en 
r\'eduction une composante de genre 2. On calcule $N:=N_7(F^{[1]}(Y)=Y^{48}+Y^{32} \mod 2$; puisque 
$(F^{[3]}(Y)=Y^4 \mod 2$; il suit que la composante de genre 1 est induite par le changement $X=y+\rho T$ o\`u $y=1 \mod 2$ est racine de  $N_7(F^{[1]}(Y)$ et $v(\rho)=(2/3)v(2)$. 

\item {\bf $\Delta_{020}$} 

On va montrer que $F(X)=1+2X^3+2^{1/2}X^5+X^7$ convient. 
Notons que  $(\bar F(X))'=X^6$ et que le changement $X=2^{2/3}T$ induit en r\'eduction une composante de genre 1; enfin le changement $X=2^{1/4}T$ induit une \'equation $Y^2=1+2^{7/4}(T^3+T^5+T^7)$ et donc un $\aalpha_2$ torseur $Y^2= T^3+T^5+T^7=f$ et puisque $df=T^2(1+T+T^2)^2dT$ on pr\'evoit que l'on a
l\`a la composante qui supporte les 3 composantes elliptiques: prouvons cela. On calcule la $7$-d\'eriv\'ee $N:=N_7(F^{[1]}(Y))$ et on analyse son polygone de Newton, quelques congruences permettent de conclure: $N =y^{48} \mod 2$, $ \rem(N,y^{33},y)=16*y^32 \mod 2^5$, $\rem(N,y^3,y)=8192*y^2 \mod 2^{14}$, $\rem(N,y^{17},y)=256*y^{16} \mod 2^9$. Il suit que les racines se r\'epartissent en 16 de valuation $\geq 5/14$ puis 
16+16 en 2 classes modulo $2^{1/4}$. 

\item {\bf $\Delta_{021}$}

On prend $F(X)=1+2^{1/2}X^5+X^7$; alors $(\bar F(X))'=X^6$. Le changement $X=2^{3/10}T$ induit en 
r\'eduction une composante de genre 2. Le changement $X=2^{1/4}T$ induit une \'equation $Y^2=1+2^{7/4}(T^5+T^7)$ et donc un $\aalpha_2$ torseur $Y^2= T^5+T^7=f$ et puisque $df=T^4(1+T)^2dT$; comme pr\'ec\'edemment on v\'erifie que l'on a ainsi la composante qui supporte les 2 composantes de genre respectif $1$ et $2$. 

\item {\bf $\Delta_{022}$}

$F(X)=1+X^7$ convient!

\end{enumerate}

\end{enumerate}
\subsubsection 
{Illustrations dans le cas $p>2$.} 
\begin{enumerate}[A.]
\item {\bf Un exemple instructif.}

Dans cet exemple $p=3$, $F(X)=1+cX^3+X^4$; pour $c\in R$. Notons que le discriminant 
$ \disc F(X)=256-27 c^4=1 \mod 3$; ainsi la g\'eom\'etrie du lieu de branchement est \'equidistante.
Dans ce cas nous avons $m-1=p$ et nous allons voir que le crit\`ere de bonne r\'eduction de 
[Le 2] est insuffisant (cf. th. 3.3.2). Nous remarquons que $F'(X)=X^2(3c +4X)$; suivant le th\'eor\`eme 3.3.2, si la courbe $Y^p=F(X)$ 
a bonne r\'eduction il existe  $d\in R$ tel que 
 $$A'(X)=(X-d)^3 \ \mod \lambda ^{3/4}$$
Comme $3 \cong \lambda^2$ il suit que la condition de [Le 2] est \'equivalente \`a 
$d^3=0 \ \mod \lambda ^{3/4}$ et donc $d =0  \mod \lambda ^{1/4}$. On doit donc tester si 
il y a bonne r\'eduction pour $ \frac{\lambda ^{3/4}}{X-d}$ et $d\lambda ^{-1/4}=0 \mod \lambda ^{1/2}$ i.e. 
$v(d-d')\geq (3/8) v(3)$. 

Appliquons l'algorithme. 

$$F(X+Y)= X^4+(c+4Y) X^3+(3 c Y+6 Y^2) X^2+(3 c Y^2+4 Y^3) X+1+c Y^3+Y^4$$
la $3$-d\'eriv\'ee de niveau $4$ vaut 

$N_4(F^{[1]}(Y))=-3^3s_0(Y)^2s_3(Y)+s_1(Y)^3 =-44 Y^9-99 c Y^8-54 c^2 Y^7-216 Y^5-270 c Y^4-54 c^2 Y^3-108 Y-27 c$

Soit  $y=d\in R'$ une racine de  $N_4(F^{[1]}(Y)$, alors $v(d)>0$ et donc 
$s_0(y)$ est une unit\'e. Alors 
$F(X+d)= X^4+(-3 s_0(y)^{1/3} s_3(y)^{2/3}+s_2(y)) X^2+(-3 s_0(y)^{2/3} s_3(y)^{1/3}+s_1(y)^3) X
+(s_0(y)^{1/3}+s_3(y)^{1/3} X)^3=
(s_0(y)^{1/3}+s_3(y)^{1/3} X)^3+(-3 s_0(y)^{1/3} s_3(y)^{2/3}+s_2(y)) X^2+ X^4.$

Notons que 
$s_2(y)=(3 c d+6 d^2)$ est divisible par 3; il suit que l'on a bonne r\'eduction pour 
$X= \lambda ^{3/4}T$.

\hfill  ///
\item
{\bf Des limites de l'algorithme.}

Dans [Ma], on montre que pour 
$a_1,a_2,...,a_n \in \ZZ_2^{nr}$ suffisament g\'en\'eraux et 
$$F(X):=\prodl{(\epsilon_1,..,\epsilon_n)\in \{0,1\}^n}
(1+(\suml{1\leq i\leq n}\epsilon_ia_i)^2X)$$ il existe $b_0,b_1,...,b_{n-1} \in \ZZ_2^{nr}$ tels que 
$$F(X)=(1+\suml{1\leq i\leq n-1}b_{i}X^{2^{n-1}-2^{i-1}})^2+b_0^2X^{2^n-1} \, \mod 4\ZZ_2^{nr}.$$
 
Cette derni\`ere congruence fournit un $2$-d\'eveloppement de Taylor de $F(X)$ \`a l'ordre 
$m=2^{n}-1$ et de niveau quelconque (cf. D\'ef. 2.4.4). 

Soit $R:=\ZZ_2^{nr}[\pi]$ avec $\pi^{2^n-1}=2^2$, alors 
l'\'equation  $Y^2=F(X)$ d\'efinit une courbe hyperelliptique ayant 
bonne r\'eduction sur $R$ relativement \`a la valuation 
de Gauss en 
$S:=(2)^{-2/(2^n-1)}X$  avec $m=2^n-1$.

Il doit \^etre clair dans l'esprit du lecteur que l'algorithme propos\'e dans cette note ne permet 
pas de montrer la congruence qui pr\'ec\`ede. 
\end{enumerate}
\section{R\'eduction stable des rev\^etements   $p$-cycliques d'une courbe de genre $>0$ qui a bonne 
r\'eduction en $p$; cas d'un lieu de branchement \`a g\'eom\'etrie \'equidistante.}
\indent{\ \ \ }
On consid\'ere un rev\^etement $p$-cyclique $C\to D:=C/G$ tel que la courbe quotient $D$ admet un mod\`ele
lisse ${\cal D}'$ et que le lieu de branchement $\Br$ se sp\'ecialise en des points distincts $\bar \Br
\subset{\cal D}'_s $. Dans ce cadre, le lemme 3.2.1 se g\'en\'eralise en:
\begin{Lem}.--- 
Apr\`es extension finie de $K$, le mod\`ele  $\cal C'$ de $C$, obtenu par normalisation du mod\`ele lisse ${\cal D}'$  de $D$ a une fibre sp\'eciale r\'eduite qui est un rev\^etement radiciel de ${\cal D}'_s$; les singularit\'es sont des cusps et se trouvent  au-dessus des z\'eros 
$\in {\cal D}'_s$ d'une forme  diff\'erentielle logarithmique $\omega$ r\'eguli\`ere  en dehors de $\bar \Br$. Le graphe d'intersection de la fibre sp\'eciale ${\cal C}_s$ du  mod\`ele stable $\cal C$ de la courbe point\'ee $(C,\Br)$  
est un arbre (il n'y a pas de cycles)  et le mod\`ele ${\cal D}$ de $D$ obtenu par quotient de l'action de $\ZZ/p\ZZ$ sur $\cal C$ a une fibre sp\'eciale qui est un arbre de 
droites projectives attach\'ees  \`a ${\cal D}'_s$ en les points de $\bar \Br\subset{\cal D}'_s $.
\end{Lem}
 
\noindent 
{\it Preuve.} Consid\'erons l'anneau local de ${\cal D}'$ au point g\'en\'erique de la fibre 
sp\'eciale; c'est un anneau de valuation discr\`ete; soit $v$ la valuation correspondante.  
Soit ${\cal C}'$, le mod\`ele  de $C$, obtenu par normalisation du mod\`ele lisse ${\cal D}'$ de $D$. 
Par [Ep], apr\`es une extension finie de $K$, on peut supposer que la fibre sp\'eciale de ${\cal C}'$ 
est r\'eduite. 

Nous allons voir qu'elle est int\`egre et que le $\mmu_p$-torseur au-dessus de $D'-\Br$ induit un 
$\mmu_p$-torseur au-dessus de ${\cal D}'_s-\bar \Br$. 

Une \'equation du rev\^etement est $Z^p=F\in K(D)$ o\`u $v(F)=0$; alors $F$ d\'efinit un diviseur principal de la $R$-courbe ${\cal D}'$, $(F):=\suml{1\leq i\leq m+1}
e_iV(\{x_i\})+p(D_0)$ o\`u $\Br =\{x_i,1\leq i\leq m+1\}$ et $V(\{x_i\})$ d\'esigne la fermeture sch\'ematique de $\{x_i\}$, enfin   $e_i\in \lv N$ avec $(e_i,p)=1$ et $(D_0)$ est 
un diviseur horizontal. On a donc modulo $\pi$; $(\bar F)=\suml{1\leq i\leq m+1}e_i\bar x_i
+p({D_0}_s)$ (cf. [Li] Lemma 7.1.29); ainsi $\bar F$ n'est pas une puissance $p$-i\`eme dans 
$k({\cal D}'_s)$ et donc l'\'equation $\bar Z^p=\bar F\in k({\cal D}'_s)$ d\'efinit un 
$\mmu_p$-torseur au-dessus de ${\cal D}'_s-\bar \Br$. La forme diff\'erentielle logarithmique 
$\omega:=d\bar F/\bar F$ est une forme diff\'erentielle sur ${\cal D}'_s$ (ind\'ependante de $F$) qui est r\'eguli\`ere en dehors de $\bar \Br$; c'est la forme diff\'erentielle associ\'ee au torseur.

La fibre sp\'eciale ${\cal C}'_s$ est alors hom\'eomorphe \`a ${\cal D}'_s$ et puisque le diviseur 
de $\bar F$ a une multiplicit\'e en $\bar x_1$ premi\`ere \`a $p$; une modification de 
l'\'equation du torseur va donner la multiplicit\'e $1$ et induire un mod\`ele lisse en $\bar x_1$. 
Ainsi les singularit\'es de la courbe ${\cal C}'_s$ sont des cusps situ\'es dans les 
z\'eros de $\omega$ et en dehors des $\bar x_i$. 

Si $x$ est  un z\'ero de $\omega$, apr\`es localisation en $x$ on se retrouve dans la situation 
locale d\'ecrite par Raynaud dans ([Ra1], d\'emonstration du th\'eor\`eme 1 p. 182)), le lemme suit. 
\hfill  ///

Le th\'eor\`eme 3.3.1 a une g\'en\'eralisation que nous d\'ecrivons bri\`evement.

Soit $\bar d\in {\cal D}'_s $, un z\'ero de $\omega$ d'ordre $m-1=m(\bar d)-1$. Une \'equation du 
torseur 
au-dessus de la fibre formelle en $\bar d$ est $Z^p=F(X)\in R[[X]]$;  o\`u $F(X)$ est une unit\'e 
modulo $\pi$. Dans le s\'epar\'e compl\'et\'e du module de diff\'erentielle $\Omega_{R[[X]]/R}$ on 
a alors l'\'egalit\'e $ d\bar F/\bar F=u(X)X^{m-1}dX\mod \pi$ o\`u  $u(X)$ est une  unit\'e modulo $\pi$ . Comme dans le cas des polyn\^omes, on fait un $p$-d\'eveloppement de Taylor sp\'ecial
\`a l'ordre $m$ et de niveau $n$ avec $p^{n}<m<p^{n+1}$ (cf. 2.5.1),

$$ F(X+Y)=E(X,Y)^p+\suml{1\leq j\leq m,\ (j,p)=1}F^{[j]}(Y)X^j \mod (p^{c_n}X,X^{m+1})$$
o\`u $E(X,Y)\in (R[[Y]])^{alg}[X]$,  $\deg_X E(X,Y)=[m/p]$ et $F^{[j]}(Y)\in (R[[Y]])^{alg}$ (notez que la congruence impose l'\'egalit\'e $E(0,Y)^p=F(Y)$). De mani\`ere analogue \`a III.4.1, on d\'efinit 
la $p$-d\'eriv\'ee de niveau $m$, $$N_m(F^{[1]}(Y)):=N_{L_n/K(Y)}(F^{[1]}(Y))\in R[[Y]].$$ 
Par le th\'eor\`eme de pr\'eparation de Weierstrass $N_m(F^{[1]}(Y))=P_m(Y)U_m(Y)$ o\`u $P_m(Y)\in R[Y]$ est un polyn\^ome distingu\'e, $P_m(Y)=Y^{(m-1)p^{r(m)}} \mod \pi$, et $U_m(Y)\in R[[Y]]$ est une unit\'e. Les
z\'eros de $P_m(Y)$ donnent des centres des disques ferm\'es de la droite projective qui induisent 
une composante de genre non nul dans la fibre sp\'eciale ${\cal C}_s$ du  mod\`ele stable $\cal C$ de la courbe point\'ee $(C,\Br)$ au-dessus du z\'ero  $\bar d\in {\cal D}_s$ de $\omega$.

Si $D$ n'est pas 
la droite projective, il y a plusieurs obstacles pour rendre les arguments pr\'ec\'edents algorithmiques.


\begin{flushleft}

Michel MATIGNON \\
Laboratoire de Th\'eorie des Nombres\\
et d'Algorithmique Arithm\'etique\\
UMR 5465 CNRS \\
Universit\'e de Bordeaux I \\
351 cours de la Lib\'eration, \\
33405 Talence Cedex, France \\
e-mail : {\tt matignon@math.u-bordeaux.fr}

\end{flushleft}


\begin{thebibliography}{Gr-Ma}

\bibitem[Ab]{Ab}
A.~Abbes, \emph{R\'eduction semi-stable des courbes d'apr\`es Artin, Deligne, Grothendieck, Mumford, Saito, Winters, $\ldots$}, Courbes semi-stables et groupe fondamental en g\'eom\'etrie alg\'ebrique (Luminy, 1998), Progr. Math., {\bf 187} (2000) Birkh\"auser,  59--110.
  
\bibitem[Co]{Co} 
R.~F.~Coleman, \emph{Computing stable reductions},  S\'eminaire de Th\'eorie des Nombres, Paris
 1985--86, Progress in Mathematics {\bf 71} (1987), Birkha\"user,  1--18.  

\bibitem[Ep]{Ep} 
H.~P.~Epp, \emph {Eliminating wild ramification}, Invent. Math.  {\bf 19} (1973), 235-249.

\bibitem[Gr-Ma ]{Gr-Ma } 
B.~Green, M.~Matignon, \emph{Order $p$ automorphisms of the open disc of a $p$-adic field}, J. Amer. Math. Soc. 12 (1999), pp. 269-303.

\bibitem[He]{He} Y.~Henrio, \emph{Arbres de Hurwitz et automorphismes d'ordre p des disques et des couronnes p-adiques formels },  Compositio Mathematica, to appear.


\bibitem[Le 1]{Le1} C.~Lehr, \emph{ Reduction of wildly ramified covers of curves}, Th\`ese University of
Pennsylvania, (2001).

\bibitem[Le 2]{Le2} C.~Lehr, \emph{ Reduction of $p$-cyclic Covers of the Projective Line},  Manuscripta Math. 106 (2001) 2, 151-175.

\bibitem[Li 1]{Li1} Q.~Liu, \emph{ Algebraic Geometry and Arithmetic Curves}, Oxford Graduate Texts in 
Mathematics, 6 (2002), Oxford University Press.

\bibitem[Li 2]{Li2} Q.~Liu, \emph{Courbes stables de genre 2 et leur sch\'ema de modules}, 
Mathematische  Annalen 295  (1993),  201-222.
  
\bibitem[Ma]{Ma} M.~Matignon, \emph{$p$-groupes ab\'eliens de type $(p,\cdots,p)$ et disques ouverts $p$-adiques}, Manuscripta Math. 99 (1999), no. 1, 93--109.

\bibitem[Mi]{Mi} J.~S.~Milne, \emph{ \'Etale Cohomology}, Princeton Mathematical
Series, {\bf 33}, (1980).

\bibitem[Ra 1]{Ra 1} M.~Raynaud, \emph{$p$-groupes et r\'eduction semi-stable des courbes}, The 
Grothendieck Festschrift, Vol III, Progress in Mathematics {\bf 88} (1990), Birkha\"user, 179--197.


\bibitem[Ra 2]{Ra 2}M.~Raynaud, {Mauvaise r\'eduction des courbes et $p$-rang}, C.R. Acad. Sci. Paris,  316, S\'erie I, (1994), 1279--1282.

\bibitem[Sa 1]{Sa 1} M.~Sa\"\i di, \emph{Galois covers of degree p: semi-stable reduction and Galois action}, math.AG/0106249.

\bibitem[Sa 2]{Sa 2} M.~Sa\"\i di, \emph{Wild ramification and a vanishing cycles formula},
math.AG/0106248.

\bibitem[Sa 3]{Sa 3} M.~Sa\"\i di, \emph{Torsors under finite and flat group schemes of rank p with Galois action}, math.AG/0106246.
\end{thebibliography}
\end{document}